\pdfoutput=1

\documentclass[final]{siamltex} 

\usepackage{graphicx}
\usepackage{amsmath}
\usepackage{amsfonts}
\usepackage{amssymb}
\usepackage{color} 
\usepackage{url} 
\usepackage[boxed]{algorithm}
\usepackage{algorithmicx}
\usepackage{algpseudocode}
\usepackage{multirow}

\newcommand{\eq}[1]{\begin{equation}\label{#1}}
\newcommand{\en}{\end{equation}}
\def\nref#1{(\ref{#1})}

\def\intv{$[\xi, \ \eta]$}

\def\trans{^{\mathsf{T}}}

\def\invt{^{-\mathsf{T}}}%
\providecommand{\norm}[1]{\Vert #1 \Vert} 
\def\inv{^{-1}}%
\def\inv{^{-1}}

\def\RR{\mathbb{R}}

\usepackage{listings}
\definecolor{hellgelb}{rgb}{1,1,0.85}
\definecolor{colKeys}{rgb}{0,0,1}
\definecolor{colIdentifier}{rgb}{0,0,0}
\definecolor{colComments}{rgb}{1,0,0}
\definecolor{colString}{rgb}{0,0.5,0}

\begin{document}

\title{The Eigenvalues Slicing Library (EVSL): Algorithms, 
Implementation, and Software}

\author{
Ruipeng Li\thanks{Center  for Applied  Scientific Computing,
    Lawrence  Livermore National  Laboratory,  P. O.  Box 808,  L-561,
    Livermore,   CA  94551   ({\tt li50@llnl.gov}).  This   work  was
    performed under the  auspices of the U.S. Department  of Energy by
    Lawrence    Livermore   National    Laboratory   under    Contract
    DE-AC52-07NA27344 (LLNL-JRNL-746200).}
\and
Yuanzhe Xi\thanks{Department of Computer Science and Engineering, University of Minnesota Twin Cities,
Minneapolis, MN 55455
({\tt \{yxi,erlan086,saad\}@umn.edu}).
Work supported by NSF under grant CCF-1505970   
and by the Minnesota Supercomputer Institute.}
\and
Lucas Erlandson\footnotemark[2]
\and Yousef Saad\footnotemark[2]
}

\maketitle

\begin{abstract}
  This paper describes a software package called EVSL (for EigenValues
  Slicing Library)  for solving large sparse  real symmetric standard
  and  generalized eigenvalue  problems.   As its  name indicates,  the
  package  exploits  spectrum slicing,  a  strategy  that consists  of
  dividing the spectrum  into a number of  subintervals and extracting
  eigenpairs from  each subinterval independently. In  order to enable
  such a  strategy, the methods implemented  in  EVSL rely on  a quick
  calculation of the  spectral density of a given  matrix, or a matrix
  pair.   What   distinguishes  EVSL  from  other    currently
  available packages is that  EVSL  relies entirely  on filtering  techniques.
 Polynomial  and rational  filtering   are both implemented
  and   are coupled  with Krylov subspace  methods and the
  subspace iteration algorithm.  On the implementation side, the package
  offers  interfaces for  various  scenarios  including   matrix-free
  modes, whereby the  user can supply his/her own  functions to perform
  matrix-vector  operations or  to solve  sparse linear  systems.  The
  paper describes  the algorithms in  EVSL, provides  details on
  their implementations, and discusses performance issues 
  for the various methods.
\end{abstract}


\section{Introduction\label{sec:intro}}
A good number of software packages have been developed in the past two
decades for  solving large  sparse symmetric eigenvalue  problems see,
e.g.,                            \cite{Baglama:2003:AIM:838250.838257,
  Baker:2009:ASN:1527286.1527287,              JADAMILU,Filtlan-paper,
  SLEPc,Stathopoulos10,lsy:ARPACK98}. In most cases,
 these packages deal with the situation where a  relatively
small number of eigenvalues are sought on  either end of the spectrum.  Computing
eigenvalues located  well inside  the spectrum  is often  supported by
providing  options in the software that  enable one to exploit spectral  transformations,
i.e., shift-and-invert  strategies~\cite{Parlett-book}.  
However,  these packages are  generally not designed  for handling
the situation when a large number  of these eigenvalues, reaching in the
thousands or tens of thousands, are sought and when, in addition,
they are   located at the interior of the spectrum.
Yet, it is now this breed
of   difficult  problems   that   is  often   encountered  in   modern
applications.  For  example, in  electronic structure  calculations, a
method such as Density Functional Theory (DFT) will see the number of ground
states increase to very large numbers in realistic simulations.  As an
illustration,     the    2011     Gordon     Bell    winning     paper
\cite{HasegawaGordonBell2011} showed  a calculation on  the K-computer
that  had as  many  as 107,292  atoms leading  to  the computation  of
229,824  orbitals (eigenvectors  of the  Kohn-Sham equation).   In the
introduction of the  same paper the authors pointed  out that \emph{``in
  order to represent  actual behavior of genuine  materials, much more
  computational resources, i.e., more CPU  cycles and a larger storage
  volume, are  needed to make  simulations of up to  100,000 atoms.''}
As a side note, it is  worth mentioning that the article discusses how
the orbitals  were divided  in 3  sets of  76,608 orbitals  handled in
parallel,  a perfect  example  of spectrum  slicing.  In  computations
related  to  excited states,  e.g.,  with  the Time-dependent  Density
Functional   Theory (TDFT)  approach,   the  situation   gets   much   worse
\cite{Russ-tdlda} since one needs to  compute eigenvalues that are not
only  related to  occupied states  but  also those  associated with  a
sizable number of  unoccupied ones.  EVSL is  specifically designed to
tackle these  challenges and  to address   standard and
generalized symmetric eigenvalue problems, encountered in large
scale applications of the type just discussed.

 As a background we will
 begin with  a brief review of the main  software packages currently
available  for large  eigenvalue problems,  noting that it is beyond the
scope of this paper to provide an exhaustive survey.
We mentioned earlier that  most of the
packages that are available today 
aim  at extracting a  few eigenvalues.  Among  these, 
the best known is undoubtedly 
ARPACK which  has now become a de-facto standard for
solvers  that  rely  on  matrix-vector operations.   ARPACK  uses  the
Implicitly Restarted  Arnoldi Method  for the non-symmetric  case, and
Implicitly   Restarted  Lanczos   Method   for   the  symmetric   case
\cite{lsy:ARPACK98}.   ARPACK  supports  various  types  of  matrices,
including single  and double  precision real arithmetic  for symmetric
and  non-symmetric  cases, and  single  and  double precision  complex
arithmetic for standard and generalized problems.  In addition, ARPACK
can also   calculate  the Singular Value Decomposition (SVD) of a matrix.  
Due to  the  widespread use  of
ARPACK, once development  stalled, vendors began creating their own
additions to  ARPACK, fracturing what used to be a uniform standard.  
ARPACK-NG  was created  as  a joint  project
between Debian, Octave and Scilab, in an effort to contain  this split.

To cope with the large memory  requirements required by the use of the
(standard)  Lanczos, Arnoldi  and Davidson  methods, restart  versions
were developed  that allowed to  reuse previous information  from the
Krylov subspace in clever  ways.  Thus, TRLan \cite{TRLan,trlan99} was
developed as a  restarting variant of the Lanczos  algorithm.  It uses
the thick-restart Lanczos method \cite{stsw:98}, whose generic version
is  mathematically  equivalent  to the  implicitly  restarted  Lanczos (IRL)
method. The TRLan    implementation  supports  
both   single  address  memory   space  and
distributed      computations
\cite{Yamazaki:2010:APS:1824801.1824805,wusi:00} and also   
supports  user provided matrix-vector product
functions. More recently, TRLan was
rewritten in  C as  the $\alpha$-TRLAN package,  which has  an added
feature    of   adaptively    selecting   the    subspace   dimensions
\cite{Yamazaki:2010:APS:1824801.1824805}.

As part  of the Trilinos project~\cite{Trilinos},  Anasazi provides an
abstract framework to solve eigenvalue  problems, enabling the user to
choose        or       provide        their       own        functions
\cite{Baker:2009:ASN:1527286.1527287}.  It utilizes advanced paradigms
such as object-oriented programming and  polymorphism in C++, to 
enable  a robust  and  easily  expandable framework.   With  the help  of
abstract interfaces,  it lets users  combine and create software  in a
modular fashion, allowing them to construct a precise setup needed for
their specific  problem. A notable  recent addition to Anasazi  is the
TraceMin
eigensolver~\cite{SamehTong00,TraceMinPar,sw:tracemin82}, which  can be
viewed     as     a     predecessor     of     the     Jacobi-Davidson
algorithm~\cite{Jacobi-Davidson}.

The BLOPEX  package \cite{Blopex}  is available as  part of  the Hypre 
library \cite{hypre},  and also  as  an  external package  to  PETSc \cite{knyazev_block_2007}. Through  this
integration,  BLOPEX  is  able  to make  use  of  pre-existing  highly
optimized  preconditioners.   BLOPEX  uses  the  Locally  Optimal  Block
Preconditioned Conjugate Gradient method (LOBPCG) \cite{LOBPCG}, which
can be combined  with a user-provided preconditioner.  It  is aimed at
computing  a few  
extreme  eigenvalues of  symmetric
positive   generalized   eigenproblems,   and   supports   distributed
computations  via the message
passing interface (MPI).  

SLEPc~\cite{SLEPc}  is a  package aimed  at providing  users with  the
ability to  use a variety of eigensolvers to solve
different types of eigenvalue problems. 
SLEPc, which is integrated into PETSc,   supports quite a few 
eigensolvers, with the ability to solve  SVD problems, 
Polynomial  Eigenvalue problems,  standard and generalized problems,
and it also  allows  solving interior eigenvalue problems with 
shift-and-invert transformations. 
SLEPc  is written in a mixture of
C,  C++,  and  Fortran,  and can be run in   parallel    on  many
different platforms.

PRIMME~\cite{Stathopoulos10} is a relatively recent
software package that aims to  be a
robust, flexible and efficient,  `no-shortcuts' eigensolver.  It includes
a multilevel interface,  that allows both experts  and novice users
to  feel at ease.  PRIMME  uses nearly  optimal restarting  via thick
restart  and  `+k restart'  to  minimize  the  cost of  restarting.   By
utilizing  dynamic switching  between methods  it is  able to  achieve
accurate results very quickly \cite{Stathopoulos10}. 
The interesting article \cite{VomelAl08} shows a (limited) comparison 
of a few  software packages and found the GD+K algorithm from PRIMME to be
the best in terms of speed and robustness 
for the problems they considered. These problems originate from 
quantum dot and quantum wire simulations in which the potential
is given and, to paraphrase the authors, 
\emph{`restricts the computation to only a small number of
interior  eigenstates  from  which  optical  and  electronic
properties can be  determined'}. In fact, the goal is to compute 
the two eigenvalues that define valence and conduction bands.
There may be thousands of eigenvalues to the left of these two eigenvalues
and the paper illustrates how  some codes, e.g., the 
IRL 
routine from ARPACK faces severe difficulties in this case
because they are designed for computing
smallest (or largest) eigenvalues. 

While the  packages mentioned  above are used  to calculate  a limited
number   of   eigenvalues,   the   FILTLAN~\cite{Filtlan-paper},   and
FEAST~\cite{FEAST}  packages are,  like  EVSL, designed  to compute  a
large number of eigenpairs associated with eigenvalues not necessarily
located  at  the extremity  of  the  spectrum.   Rather than resort to 
traditional shift-invert methods, FILTLAN uses a combination of
(non-restarted) Lanczos
and polynomial  filtering with Krylov projection  methods to calculate
both  the  interior and  extreme  eigenvalues.   In addition, it  uses
the Lanczos algorithm  and  partial  reorthogonalization to improve  performance
\cite{Filtlan-paper}.  
FEAST has been designed with the same motivation as EVSL, namely to solve
the kind of eigenvalue problems  that are prevalent in solid state physics.
It can  search for an arbitrary number
of eigenvalues of Hermitian and non-Hermitian  eigenvalue problems.
Written in Fortran, it exploits the well-known 
Cauchy integral formula to express
the eigenprojector. This formula is approximated via numerical integration
and this leads to a rational filter to which subspace iteration is applied.
Single  node and  multi-node versions are  available via
OpenMP and MPI.
The idea of using Cauchy integral formulas has also been exploited by 
Sakurai and co-authors~\cite{SakSug03,SakTad07,Asakura2009}.
A related software package called z-Pares 
\cite{zPares} has been made  available  from the University of Tsukuba,
which is implemented in Fortran 90/95 and
 supports single and double precision. It offers interfaces for 
both sparse and dense matrices but also allows arbitrary matrix-vector
products via reverse communication. This code was developed with
a high degree of parallelism in mind and, for example, it can utilize
a 2-level distributed parallelism via a pair of MPI communicators.

A number of other packages implement specific classes of methods.
Among these is 
JADAMILU \cite{JADAMILU} which focuses on the Jacobi-Davidson framework
to which it adds a battery of preconditioners. 
It too is geared toward the  computation of a few selected eigenvalues.
The above list is by no means exhaustive and the field is constantly
evolving. What is also interesting to note is that these packages are 
rarely used in applications such as the ones mentioned above
where a very large number of eigenvalues are computed. 
Instead, most application packages, including 
PARSEC \cite{PARSEC-paper}, CASTEP \cite{CASTEP},
ABINIT \cite{ABINIT}, Quantum Expresso \cite{QE}, OCTOPUS \cite{OCTOPUS}, 
and VASP \cite{VASP2} implement their own eigensolvers, or other 
optimization schemes that bypass eigenvalue problems altogether.
For example, PARSEC implements a form of nonlinear subspace iteration
accelerated by Chebyshev polynomials. In light of this it may be argued that
efforts to develop general purpose software to tackle this class of 
problems may be futile. The counter-argument is that 
even if a small subset from a package, or a specific technique from it,
 is  \emph{adapted} or \emph{retrofitted} into the application, instead of 
the package being entirely
\emph{adopted}, the effort will be worthwhile if it leads to significant
gains.

We end this section by introducing our notation and the terminology 
used throughout the paper.
The generalized symmetric eigenvalue problem (GSEP) considered is of the form:
\begin{equation}
\label{eq:AB}
Ax =\lambda Bx,
\end{equation}
where both $A$ and $B$ are $n$-by-$n$ large and sparse real symmetric and $B$ is positive definite.  When $B=I$, \eqref{eq:AB} reduces to the 
standard real symmetric eigenvalue problem (SEP)
\begin{equation}
\label{eq:SEP}
Ax =\lambda x.
\end{equation}
For generalized eigenvalue problems we will often refer to the 
Cholesky factor $L$ of $B$ which satisfies $B = L L^T$, where
$L$ is lower triangular. We will also refer to 
$B^{1/2}$ the matrix square root of $B$.
 It is often the case that $L$ is expensive to compute, e.g.,
for problems that arise from 3-D simulations. In this cases, we may
utilize  the factorization $B = B^{1/2} \cdot B^{1/2}$ and approximate 
the action of  $ B^{1/2}$ by that of a low degree polynomial in $B$.
If the smallest eigenvalue is $\lambda_{min} $ and the largest
$\lambda_{max}$, we refer to the interval $[\lambda_{min}, \ \lambda_{max} ]$
as the \emph{spectral interval}. We will often use this term for 
an interval $[a, \ b]$ that tightly includes 
$[\lambda_{min}, \ \lambda_{max} ]$.

\section{Methodology}
This  section begins with an  outline of the methodologies of  the spectrum
slicing algorithms  as well as the projection methods   for eigenvalue problems
that have been  implemented in EVSL. For further reading, the theoretical 
and algorithmic
details
can be found in \cite{spectrumslicing,gdos,lsfeast}.

\subsection{Spectrum slicing}
\begin{figure}
\caption{An illustration of slicing a spectrum into $5$ subintervals $[t_i, t_{i+1}]$ $(i=0,\ldots, 5)$.  The solid blue curve represents a smoothed DOS and the dotted red lines separate the subintervals.}
\centerline{\includegraphics[width=0.62\textwidth]{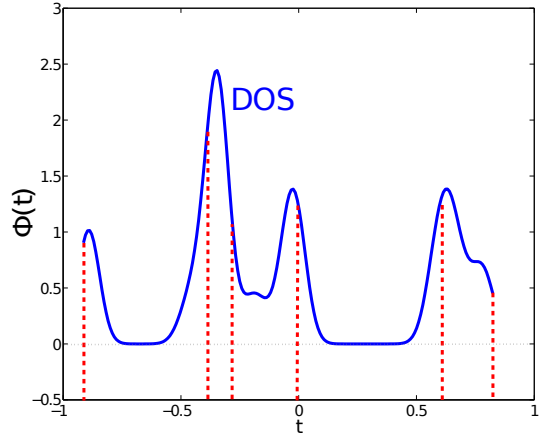}}
\label{fig:dosslice}
\end{figure}

The EVSL package relies on  a \emph{spectrum slicing} approach to deal
with  the  computation of a very  large number  of eigenvalues.  The
idea of spectrum  slicing consists of 
a \textit{divide-and-conquer} strategy:
the  target interval $[\xi, \eta]$ that contains  many eigenvalues
is subdivided into a number of 
subintervals  and  the  eigenvalues  from each  subinterval  are  then
computed  independently from  the  others.   The  simplest
slicing method would be to  divide $[\xi, \eta]$  uniformly into  equal size
subintervals.  However, when the distribution of eigenvalues is highly
nonuniform, some  subintervals may contain many more eigenvalues than
others and this  will cause a severe  imbalance in the computational
time and, especially, in  memory usage.   EVSL
  adopts a  more sophisticated  slicing strategy  based on  
exploiting the \emph{spectral density} that is also known
\emph{density   of   states}   (DOS).  
 The spectral density is a
function $\phi(t)$  of  real number $t$ that provides for any given $t$ 
a probability measure for finding eigenvalues near $t$. 
Specifically,
$n\phi(t) \Delta t$ will give an approximate number of eigenvalues 
in an infinitesimal interval surrounding $t$ and of width $\Delta t$.
For   details,   see,   e.g.,
\cite{LinYangSaad2013}. Formally,
 the DOS for a Hermitian matrix $A$ is defined 
as follows:
 \eq{eq:DOS0}   \phi(t)  =  \frac{1}{n}  \sum_{j=1}^n   \delta(t  -
\lambda_j) , \en where $\delta$  is the Dirac $\delta$-function or the
Dirac distribution, and  the $\lambda_j$'s are the  eigenvalues of $A$
or $B\inv A$.  Written in this formal form, the DOS is not a proper
function but a distribution and it is its approximation by a smooth
function that is of interest~\cite{LinYangSaad2013}.
Notice that if the DOS function $\phi$ were available,
the exact  number of eigenvalues located 
inside  $[\xi, \eta]$ could be obtained from the integral
\[ 
 \int_{\xi}^{\eta} n\phi(t) dt = \int_{\xi}^{\eta} \sum_j\delta(t-\lambda_j) dt.
\] 
Thus, the task of slicing $[\xi,\eta]$ into $n_s$ subintervals 
that roughly contain  the same number of eigenvalues can
 be accomplished by ensuring that the integral of $\phi$  on each subinterval
is an equal share of the integral on the whole interval, i.e,
\begin{equation} \label{eq:slice-ti} 
 \int_{t_i}^{t_{i+1}} \phi(t) dt = 
\frac{1}{n_{s}}\int_{\xi}^{\eta} \phi(t) dt , \quad i=0,1,\ldots, n_s-1,
\end{equation}
where
$[t_i, t_{i+1}]$ denotes each subinterval 
with $t_0 = \xi$ and $t_{n_s} = \eta$.
The values of the endpoints $\{t_i\}$ can be computed easily 
by first finely discretizing $[\xi,\eta]$ with  $N \gg n_s$      evenly       spaced  points $x_0=\xi\!<\!x_1\!<\!\ldots\!<\!x_{N-2}\!<\!x_{N-1}=\eta$ and then
placing each $t_i$ at some $x_j$ such that  \eqref{eq:slice-ti}
can be fulfilled approximately.
Figure \ref{fig:dosslice} shows  an illustration of spectrum slicing with DOS.

Clearly, the \emph{exact  DOS} function $\phi$ is  not available since
it   requires   the   knowledge   of   the   eigenvalues   which   are
unknown.  However, there  are inexpensive  procedures that  provide an
approximate  DOS   function  $\tilde  \phi$,  with   which  reasonably
well-balanced slices  can usually be  obtained.  Two types  of methods
have been developed for computing  $\tilde \phi$.  The first method is
based   on   the   kernel   polynomial  method   (KPM),   see,   e.g.,
\cite{LinYangSaad2013} and the references therein, and the  second one
uses  a   Lanczos  procedure  and  the   related  Gaussian  quadrature
\cite{gdos}.   These  two  methods   are  in  general  computationally
inexpensive  relative to  computing  the eigenvalues  and require  the
matrix $A$ only through matrix-vector operations.

The above  methods for  computing DOS can  be extended  to generalized
eigenvalue   problems,   $Ax=\lambda   Bx$.   In   addition   to   the
matrix-vector products with $A$, these methods require solving systems
with $B$, and  also with $B^{1/2}$ or with $L\trans$  where $L$ is the
Cholesky factor of $B$, to generate a correct starting vector.  A
canonical way to do this is to obtain the Cholesky factor via  a sparse
direct method. This may sometimes be too expensive 
whereas in  many applications  that arise  from finite
element discretization, where $B$ is a  very well conditioned
mass matrix, 
using  iterative   methods  such   as  Chebyshev  iterations   or  the
least-squares polynomial approximation discussed in \cite{gdos} can become
a much more attractive alternative.

\subsection{Filtering techniques}\label{sec:filters} 
\begin{figure}[tbh]
\caption{An illustration  of filtering  techniques for  real symmetric
  eigenvalue  problems.  Left:  The eigenvalues  inside  $[-1,1]$  are
  mapped to the eigenvalues of the  filtered matrix, which are the red
  dotted points along  the blue solid curve.  Right:  The spectrum of the
  filtered  matrix,  where the  wanted  eigenvalues  in $[0,0.3]$  are
  mapped into the largest ones  of the filtered matrix. In this case the
 filter has been built so that these transformed eigenvalues are 
  larger than $0.8$.}
\centerline{
\includegraphics[width=0.55\textwidth]{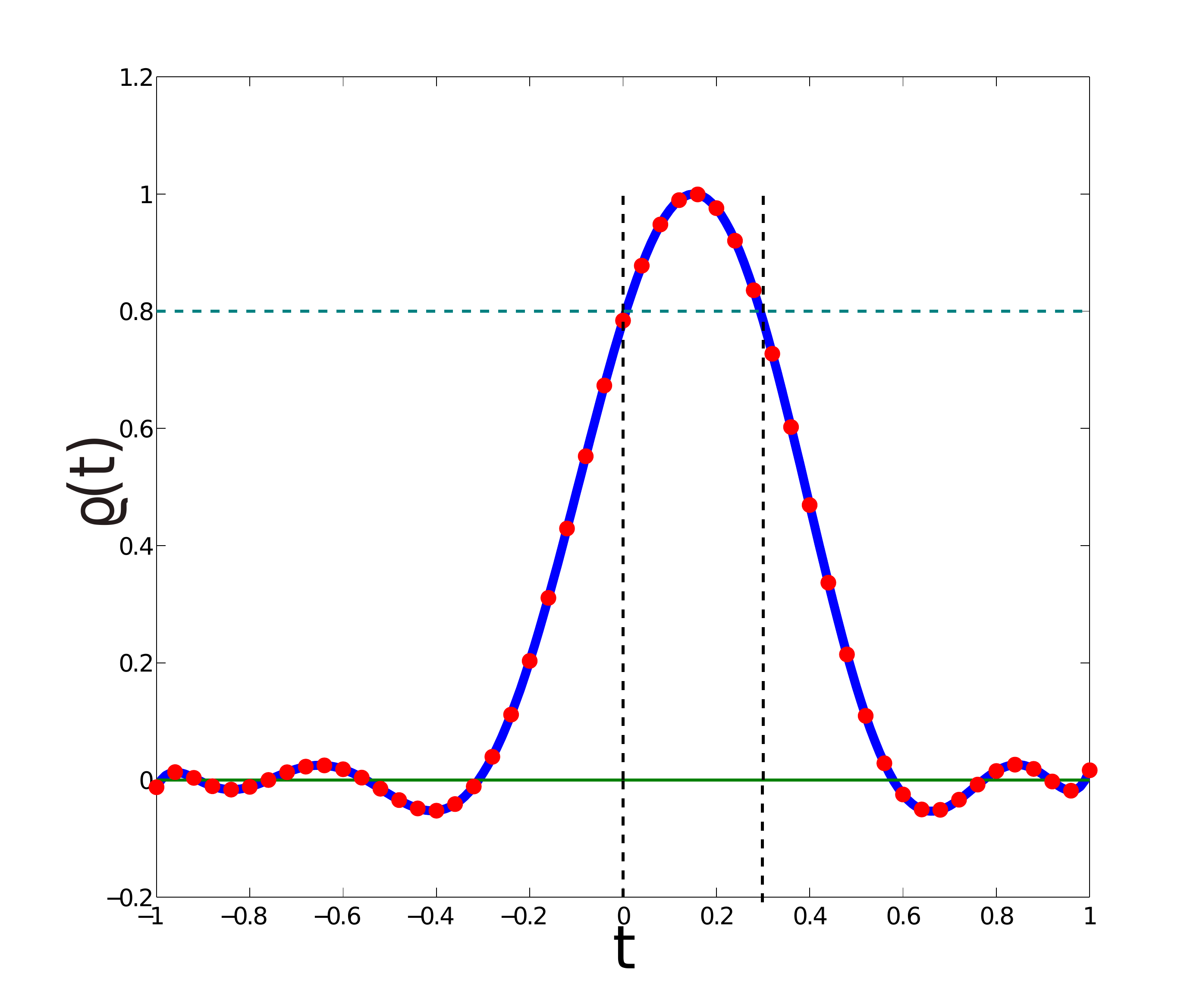}\hspace*{-6mm}
\includegraphics[width=0.55\textwidth]{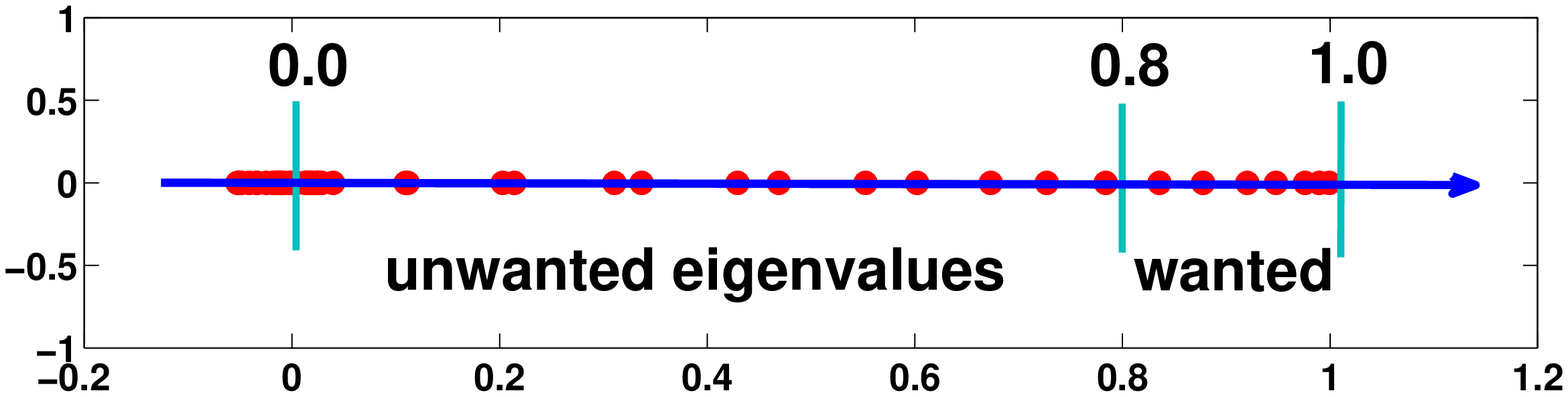}  }
\label{fig:filtering}
\end{figure}

An effective strategy for  extracting  eigenvalues 
from an interval $[t_i, t_{i+1}]$
that is  deep inside the  spectrum, is to 
apply  a  projection method  on  a  filtered  matrix, where  a  filter
function $\rho(t)$ aims at  amplifying the desired eigenvalues located
inside $[t_i, t_{i+1}]$ and dampening those  outside $[t_i, t_{i+1}]$. This is also referred to as spectral transformation and
it can be achieved via a polynomial filter or a rational one.
Both types  of filtering techniques
have  been developed and  implemented in   EVSL, with a goal of
providing different effective options for 
 different scenarios.  The main goal  of a filter function
is to  map the wanted part of  the spectrum of the
original  matrix to  the largest  eigenvalues of  the filtered  matrix.
Figure
\ref{fig:filtering}  illustrates a  polynomial filter  for the interval
$[0,0.3]$.

\subsubsection{Polynomial filtering}\label{sec:polyfit} 
The  polynomial  filter  adopted  in EVSL  is  constructed  through  a
Chebyshev polynomial  expansion of the  Dirac delta function  based 
at some well-selected point inside the
 target  interval.  Since  Chebyshev polynomials  are
defined  over the  interval $[-1,  \ 1]$,  a linear  transformation is
needed to map  the eigenvalues of $B\inv A$ onto  this interval.  This
is   achieved   by   applying    a   simple   linear   transformation:
\eq{eq:mapping}   \widehat{A}    =   \frac{A   -   c    B}{d}   \qquad
\mbox{with}\quad c  = \frac{\lambda_{\max} +  \lambda_{\min}}{2} \quad
\mbox{and} \quad  d =  \frac{\lambda_{\max} -  \lambda_{\min}}{2}, 
\en
where  $\lambda_{\max}$  and  $\lambda_{\min}$  are  the  largest  and
smallest  eigenvalues of  $B\inv A$,  respectively. These  two extreme
eigenvalues can be efficiently estimated  by performing a few steps of
the  Lanczos   algorithm.   Similarly,  the  target   interval  $[t_i,
t_{i+1}]$  should  be  transformed into  $[\hat{t}_i,  \hat{t}_{i+1}]
\equiv [(t_i-c)/d, (t_{i+1}-c)/d]$.  Let $T_j(t)$ denote the Chebyshev
polynomial  of  the first  kind  of  degree  $j$.  The  $k$-th  degree
polynomial  filter   function  $\rho(t)$  for   interval  $[\hat{t}_i,
\hat{t}_{i+1}]$  centered at  $\gamma$  takes  the form  \eq{eq:Cheb0}
\rho(t)   =   \frac{\sum_{j=0}^k  \mu_j   T_j(t)}{\sum_{j=0}^k   \mu_j
  T_j(\gamma)}, \en where the expansion coefficients are
\begin{equation}
\label{eq:mu}
\mu_j =
\begin{cases}
\frac{1}{2} & \mbox{if} \  j = 0 \\ 
\cos(j \cos\inv (\gamma) ) & \mbox{otherwise}  
\end{cases} \quad .
\end{equation}
Notice that the filter function $\rho(t)$ only contains two parameters: the degree $k$ and the 
center $\gamma$, which are determined as follows.
From a lowest degree (which is $2$), we keep increasing $k$ until 
the values $\rho(\hat{t}_i)$ and $\rho(\hat{t}_{i+1})$ 
at the boundary of the interval are less than or equal to 
 a predefined threshold $\tau$, where 
by default $\tau=0.8$ and we have $\rho(\gamma)=1$ by construction.
For any degree $k$ that is attempted, the first task is to select
 $\gamma$ so that the filter is ``balanced'' such that its values at the ends of the interval are the same, i.e., 
 $\rho(\hat{t}_i) =  \rho(\hat{t}_{i+1})$. This can be achieved by applying
 Newton's method to solve the equation 
\eq{eq:Neq} 
 \rho( \hat{t}_j) - \rho (\hat{t}_{j+1}) = 0,
\en 
where the unknown is $\gamma$.
If  $\rho(\hat{t}_i)$  and  $\rho(\hat{t}_{i+1})$ do
not exceed $\tau$, the selection procedure terminates and returns 
the current polynomial of degree $k$. Otherwise $k$ is increased and the procedure will be repeated until
a maximum allowed degree is reached. There are situations that very small 
intervals require a  degree of a few thousands.

\begin{figure} 
\caption{Polynomial filters of degree 16
 with the Lanczos $\sigma$-damping and without damping for  interval 
 $[0.2, 0.4]$.}\label{fig:twoFilt}
\centerline{\includegraphics[width=0.6\textwidth]{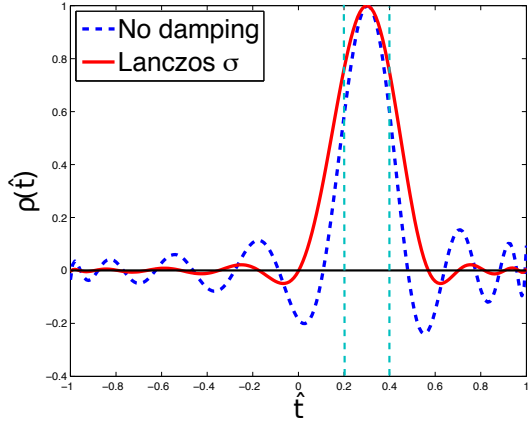}}
\end{figure}

In order to reduce the oscillations near the boundaries of
the interval, it is customary to introduce  damping multipliers to 
$\mu_j$. Several damping approaches are available in EVSL, where
the default Lanczos $\sigma$-damping~\cite[Chap. 4]{Lanczos-book} is
given by
\eq{eq:polPJ} 
\hat \mu_j = \sigma_j^{k} \mu_j,  
\en 
with
\[
\sigma_j^k  =
\begin{cases}
1 & \mbox{if} \ j=0\\
\frac{ \sin ( j \theta_k ) }{j \theta_k}, \quad
\theta_k =\frac{\pi}{k+1} & \mbox{otherwise}
\end{cases} \quad .
\]
An  illustration of  two  polynomial
filters with and without the Lanczos $\sigma$-damping
is given Figure~\ref{fig:twoFilt},  
which clearly shows the effect of damping factors for annihilating oscillations.
For  additional  details  on  selecting  the  degree, the balancing procedure,
and damping  oscillations  of the  polynomial filter, we  refer  to our
earlier paper \cite{spectrumslicing}.

\subsubsection{Rational filtering} \label{sec:ratfit}  The second type
of filtering techniques  in EVSL is rational  filtering, which can
be  a  better alternative  to  polynomial  filtering in  some  cases.
Traditional  rational  filter functions  are  usually  built from  the
Cauchy  integral  representation  of   the  step  function  $h_{[t_j,\
  t_{j+1}]}$.  Formally, the Cauchy integral of the resolvent yields a
spectral projector associated with the eigenvalues inside the interval
$[t_j,  t_{j+1}]$.  A  good approximation  of this  projector can  be
obtained  by employing  a numerical  quadrature rule.  Thus, the  step
function is approximated as follows: \eq{eq:poles} h_{[t_j, t_{j+1}]}
=   \frac{1}{2  \pi   i}   \int_{\Gamma}   \frac{1}{s-z}  ds   \approx
\sum_{j=1}^{2p} \frac{\alpha_{j}}{(z-\sigma_j)}  = \rho(z),  \en where
$\Gamma$ denotes  a circle  centered at $(t_j+t_{j+1})/2$  with radius
$(t_{j+1}-t_{j})/2$. The  poles $\sigma_j\in \mathbb{C}$ are  taken as
the quadrature nodes  along $\Gamma$ and $\alpha_j$  are the expansion
coefficients.  EVSL  implements both Cauchy integral-based  rational filters
and the   more   flexible  least-squares  (L-S)  rational  filters
discussed in
\cite{lsfeast}.   In general, a L-S rational  filter takes  the following
form:  
\eq{eq:mpoles}  \rho(z)  =  \sum_{j=1}^{2  p}  \sum_{k=1}^{k_j}
\frac{\alpha_{jk}}{(z-\sigma_j)^{k}},  \en  
where  pole $\sigma_j$  is
repeated $k_j$ times with the aim of reducing the number of poles 
on one the hand and improving the  quality of the filter for extracting
eigenvalues
 on the  other.  In particular, if $\sigma_i$  is equal to
the conjugate of $\sigma_{i+p}$,
\eqref{eq:mpoles}
can  be simplified  as \eq{eq:mpoles2}  \rho(z) =  2 \Re e \sum_{j=1}^{p}
\sum_{k=1}^{k_j}   \frac{\alpha_{jk}}{(z-\sigma_j)^k}.   \en   In  its
current implementation,  the procedure for constructing  the L-S rational
filter  starts  by choosing  the  quadrature  nodes, i.e.,  the  poles
$\sigma_j$.   Once  the  $\sigma_j$'s   are  selected,  the  coefficients
$\alpha_{jk}$  can  be  computed  by  solving  the  following  problem
\eq{eq:obj0}  \min \left\|  h_{[t_j, t_{j+1}]}  - \rho(z)\right\|_w^2
,\en  where the  norm  is  associated with  the  standard $L_2$  inner
product    \eq{eq:dotP}   \left\langle    f,    g   \right\rangle    =
\int_{-\infty}^{\infty} w(t) f(t)  \overline{ g(t) } dt ,  \en and the
weight function $w(t)$ is taken to be of the form \eq{eq:woft} w(t) =
\begin{cases}
0.01 & \mbox{if} \quad t_j\leq t\leq t_{j+1} \\
1     & \mbox{otherwise}  
\end{cases} \quad .
\en 
Since it is desirable  to have the same fixed value $1/2$ at the boundaries
 as for the Cauchy integral filters, the L-S rational filter
(\ref{eq:mpoles2}) is  scaled by $2\rho(t_j)$.

\begin{figure}[tbh]
\caption{Left: Comparison  of a  least-squares rational filter and two Cauchy integral filters using  the Gaussian Legendre rule and the  mid-point rule with $p=3$. The
  L-S rational filter uses the same poles as the Gaussian Legendre filter. 
  Right:  Comparison  of  L-S rational  filters obtained  from
  repeating one pole located at $(0,1)$ $k$ times, where $k = 2, 4, 6$.}
\centerline{
\includegraphics[width=0.54\textwidth]{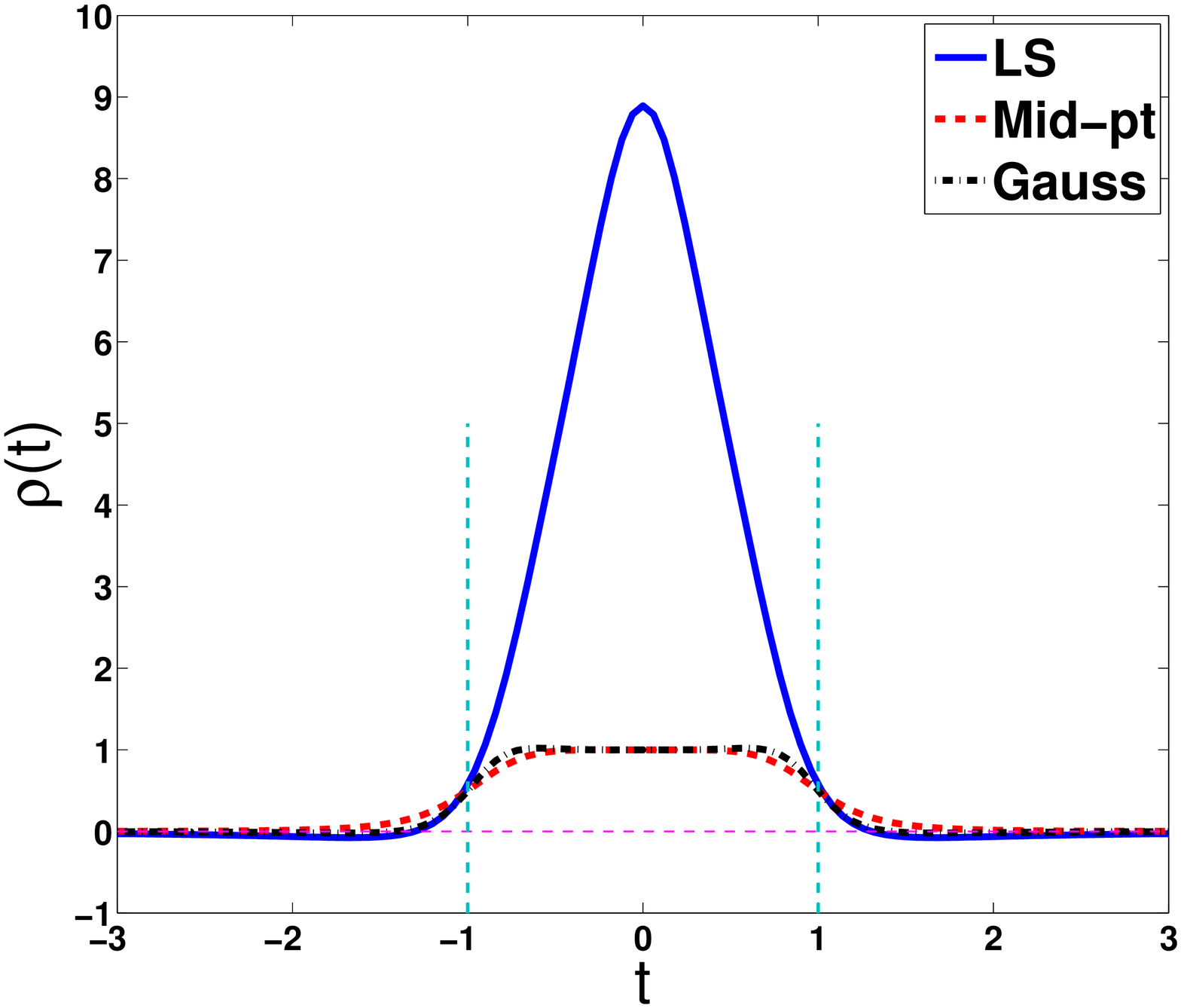}
\includegraphics[width=0.54\textwidth]{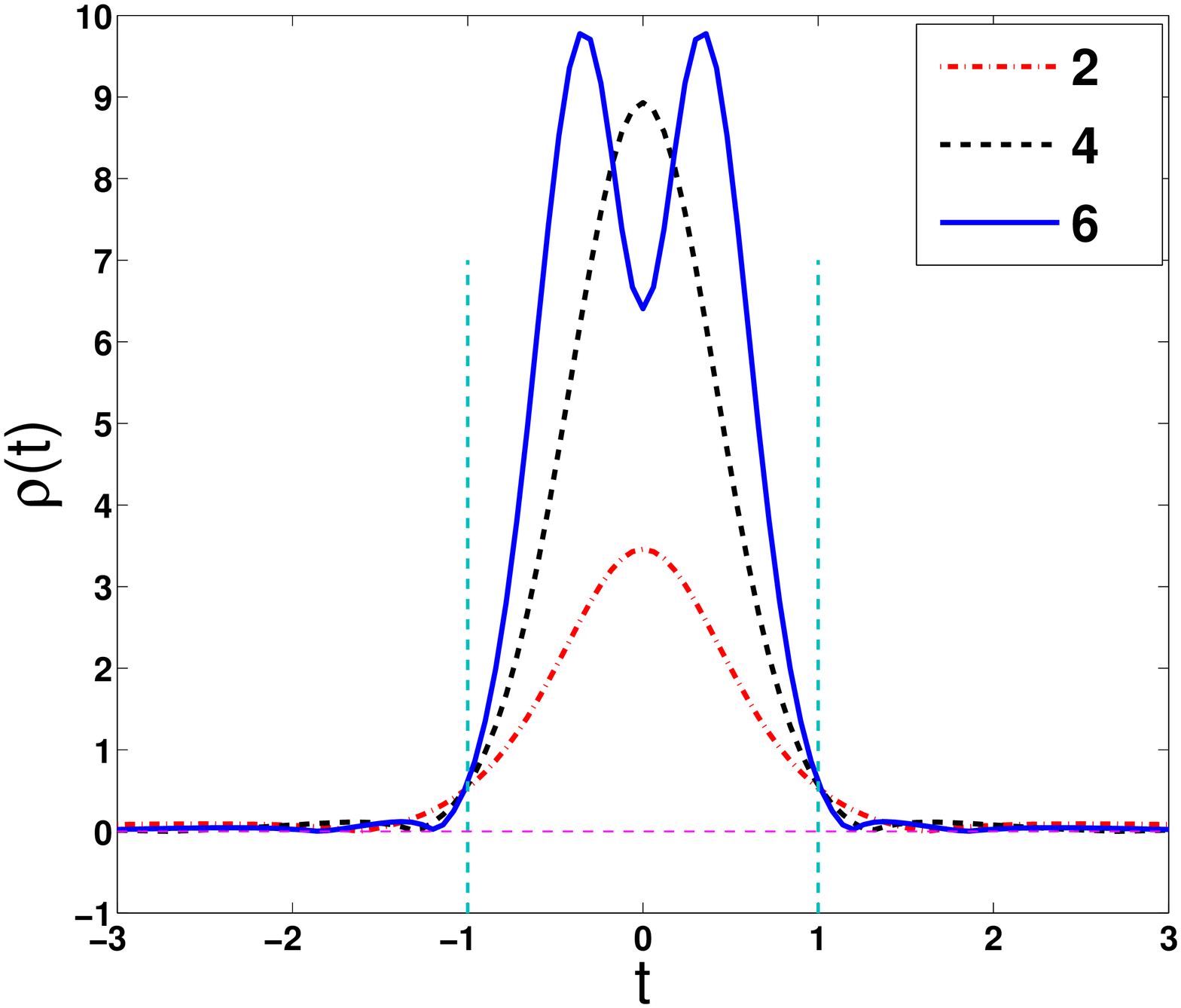} } 
\label{fig:multpole}
\end{figure}

 The left panel of Figure~\ref{fig:multpole}   shows  a 
comparison  of a L-S rational filter and two Cauchy integral  rational 
filters  for the interval $[-1, 1]$ both using  $p=3$ poles.
The Cauchy integral  rational filters are constructed  by applying the
Gauss   Legendre   rule  and   the   mid-point   rule  to   discretize
\eqref{eq:poles}, whereas  the L-S  rational filter is computed with
the  same  poles   and  the  coefficients  are   computed  by  solving
\eqref{eq:obj0}.  In this example, the  L-S rational filter  achieves almost
the same  damping effect as the  Gauss Legendre filter but  unlike the
other two rational filters, it  does not approximate the step function
and amplifies the wanted eigenvalues  to much higher values.  The
right panel shows three L-S rational filters constructed by repeating
only  $1$ pole  located at  $\sigma_1=i$ (where  $i$ is  the imaginary
unit) for  $k$ times.  It shows  that the  L-S rational  filters decay
faster across the boundaries as the value of $k$ increases.

\subsection{Lanczos algorithms for standard eigenvalue problems}
For a symmetric matrix $A$, the Lanczos algorithm builds an 
orthonormal basis of the Krylov subspace
\[
\mathcal{K}_m := \mathrm{span} \{ q_1, A q_1, \ldots, A^{m-1} q_1 \},  \]
where, at step $j$, the vector 
$Aq_j$ is  orthonormalized against $q_j$ and $q_{j-1}$ when $j>1$.
\begin{equation}
\label{eq:3term}
\beta_{i+1} q_{i+1} = A q_i - \alpha_i q_i - \beta_i q_{i-1}.
\end{equation}
In \emph{exact  arithmetic}, the 3-term recurrence would  deliver 
an orthonormal basis  of the Krylov subspace
$K_m$ but in the presence of rounding, orthogonality between the $q_i$'s 
is quickly lost, and
a form of reorthogonalization is needed in practice. 
EVSL uses full reorthogonalization to enforce the
 orthogonality  to working precision,
which is performed by no more than two steps of the classical
Gram-Schmidt algorithm  \cite{luc05,Giraud20051069} with the DGKS 
correction~\cite{10.2307/2005398}.
Let  $T_m$ denote the symmetric tridiagonal matrix
\begin{equation}
\label{eq:Tm}
T_m = \mbox{tridiag} (\beta_i, \alpha_i, \beta_{i+1}),
\end{equation}
where  $\alpha_i,\beta_i$ are from \eqref{eq:3term} and define 
$Q_m=[q_1,q_2,\ldots, q_m]$. 
Let  $(\theta_i,y_i)$ be an eigenpair of $T_m$. 
It is well-known that the
Ritz values and vectors $(\theta_i, Q_m y_i)$ can provide good
 approximations to extreme eigenvalues 
and vectors of $A$ with a fairly  small value of $m$.
EVSL adopts a simple non-restart Lanczos algorithm and a thick-restart (TR) Lanczos algorithm \cite{Andreas-al-thickr,wusi:00}
as they blend quite naturally with the filtering approaches.
In addition to TR, another essential ingredient is the inclusion of a 
``locking'' mechanism, 
which consists of an explicit deflation procedure 
for the converged Ritz vectors.

\subsection{Lanczos algorithms for generalized eigenvalue problems}
\label{sec:filtlangen} 
A straightforward way to extend the approaches developed for standard
eigenvalue problems to generalized eigenvalue problems
is to transform them into the standard form.
Suppose the Cholesky factorization $B=LL\trans$ is available, 
problem \eqref{eq:AB} can  be rewritten as
\eq{eq:stdL}
L\inv A L\invt y = \lambda y, \quad \mbox{with} \quad y=L\trans x,
\en 
which takes the standard form with the symmetric matrix $L\inv A L\invt$.
The matrix $L\inv A L\invt$ does not need to be formed explicitly,
since the procedures discussed above utilize the matrix in the form
of matrix-vector products.
However, there are situations where  
the Cholesky factorization of $B$ needed in \eqref{eq:stdL}
is not available but where one can still solve linear systems with 
the matrix $B$. In this case, we can  express 
\eqref{eq:AB} differently:
\eq{eq:stdB}
B\inv A x = \lambda x, 
\en 
which is again in the standard form but with a nonsymmetric matrix.
However, it is well-known that the matrix $B\inv A$ is self-adjoint with
respect to the $B$-inner product and this observation allows  one to use
standard methods that are designed for  self-adjoint linear operators. 
Consequently, the 3-term recurrence in the Lanczos algorithm  becomes
\begin{equation}
\beta_{m+1} w_{m+1} = B\inv A w_m -  \alpha_m w_m - \beta_m w_{m-1}. 
\end{equation}
Note here that the vectors $w_i$'s are now $B$-orthonormal, i.e., 
\[
(w_i,w_j)_B \equiv
(Bw_i,w_j) = 
\begin{cases}
0 & i \neq j \\
1 & i =j
\end{cases} \quad .
\]
The Lanczos algorithm for \eqref{eq:stdB}
is given in Algorithm~\ref{alg:LanAB} \cite[p.230]{Saad-book3},
where an auxiliary sequence $z_j \equiv B w_j$ is saved
to avoid the multiplications with $B$ in the inner products.

\begin{algorithm}[ht]
\caption{Lanczos algorithm for $Ax=\lambda Bx$\label{alg:LanAB} }
\begin{algorithmic}[1]  
\State Choose an initial vector $w_1$ with $\norm{w_1}_B=1$. Set $\beta_1=0$, $w_0=0$, $z_0=0$, 
       and compute $z_1=Bw_1$.
\For {$i=1,2,\ldots$}
\State $z:=Aw_i-\beta_i z_{i-1}$
\State$\alpha_i=(z,w_i)$
\State $z:=z-\alpha_i z_i$
\State Full reorthogonalization: $z:=z-\sum_j(z,w_j)z_j$ for $j\le i$
\State $w:=B\inv z$
\State $\beta_{i+1}=(w,z)^{1/2}$
\State $w_{i+1}:= w / \beta_{i+1}$
\State $z_{i+1}:= z / \beta_{i+1}$
\EndFor
\end{algorithmic}
\end{algorithm}

\subsection{Lanczos algorithms with polynomial and rational filtering}
Next, we  will discuss how to  combine the Lanczos algorithm  with the
filtering techniques in order to   efficiently compute all the eigenvalues
 and corresponding eigenvectors in a given  interval. 
For standard eigenvalue problem $Ax=\lambda x$,  a
filtered  Lanczos algorithm  essentially applies  the   Lanczos
algorithm with the matrix $A$  replaced by $\rho(A)$, where $\rho(A)$ 
is a filter function of $A$, in which $\rho$ is 
either a polynomial or a rational function. 
 The filtered  TR Lanczos algorithm with  locking is
sketched  in Algorithm~\ref{alg:FiltLanTR}.   In  this algorithm,  the
``candidate  set'' consists  of the  Ritz pairs  with the  Ritz values
greater  than  the filtering  threshold  $\tau$.   Then, the  Rayleigh
quotients  $\lambda_i$ with  respect to  the original  matrix $A$  are
computed.  Those  that are outside  the wanted interval  are rejected.
For the remaining candidates, we compute  the residual norm $\| A u_i-
\lambda_i  u_i \|_2$  to test  the convergence.   The converged  pairs
$(\lambda_i,  u_i)$ are  put  in  the ``locked  set''  and all  future
iterations will perform a deflation  step against those converged Ritz
vectors that are stored in the matrix  $U$, while all the others go to
the ``TR set'' for the next restart.  The Lanczos process is restarted
when the  dimension reaches  $m$.  A  test for  early restart  is also
triggered for every $Ncycle$ steps by  checking if the sum of the Ritz
values  that are  greater  than  $\tau$ no  longer  varies in  several
consecutive checks.   Once no candidate  pairs are found, we  will run
one   extra    restart in order to minimize the occurrences of missed
eigenvalues
\cite{Filtlan-paper}.  In practice, the  two stopping criteria for the
inner loop and the outer loop of Algorithm~\ref{alg:FiltLanTR} are often
found  to be effective in that they achieve the desired 
goal of neither  stopping  too early and miss a few eigenvalues,
 nor so late that  unnecessary additional work is carried out.

\begin{algorithm}
\caption{The filtered TR Lanczos algorithm with locking for 
$Ax=\lambda x$ \label{alg:FiltLanTR}}
\begin{algorithmic}[1]
\State \textbf{Input:} symmetric matrix $A \in \RR^{n \times n}$, initial unit vector $q_1$, filter function $\rho$ and interval $[t_j, t_{j+1}]$
\State \textbf{Initialization:} $q_0:=0$, $Its:=0$, $lock := 0$, $l:=0$ and $U := [\ ]$
\While{$Its \leq MaxIts$} 
\If{$l>0$}
\State Perform a special TR step  with $(I-UU\trans)\rho(A)$ and  $l:=l+1$
\EndIf
\For{$i=l+1,\dots,m$}
\State $Its :=Its+1$
\State Perform a Lanczos step with $(I-UU\trans)\rho(A)$ and full reorthogonalization
\If{$(i-l) \mod Ncycle = 0$}
\State $t_{old} = t_{new}$ and $t_{new}= \sum \theta_j$ for $\theta_j \ge \tau$
\If {$\left| t_{new} - t_{old}\right| < \tau_1$ }
\State \textbf{break}
\EndIf 
\EndIf 
\EndFor
\State Compute candidate Ritz pairs, i.e.,
$(\theta_j, u_j)$ such that $\theta_j\ge \tau$
\State Set TR set $Q := [\ ]$ and $l:=0$
\For{each candidate pair $(\theta_j, u_j)$}
\State Compute $ \lambda_j =  u_j ^T A  u_j$  
\If {$ \lambda_j \notin [t_j, t_{j+1}]$}
\State ignore this pair 
\EndIf
\If{$(\lambda_j, u_j)$ has converged}
\State Add $u_j$ to locked set $U := [U, u_j]$ and  $lock:=lock+1$
\Else
\State Add $u_j$ to TR set $Q := [Q, u_j]$ and  $l:=l+1$
\EndIf
\EndFor
\If {No candidates found in two restarts} 
\State \textbf{break}
\EndIf
\EndWhile
\end{algorithmic}
\end{algorithm}

Special formulations are needed if we are to implement 
filtering to a generalized eigenvalue problem, $Ax=\lambda Bx$.
The base problem is
\eq{eq:filtBase}
\rho(B\inv A) x = \rho(\lambda) x 
\en
which is in standard form but uses a nonsymmetric matrix.
Multiplying both sides by $B$ yields the form:
\begin{equation} \label{eq:gen1}
K x = \rho(\lambda) B x, \quad \mbox{with} \quad K=B\rho(B\inv A)
\end{equation}
where we now have a generalized problem with a symmetric matrix $K$ and
an SPD matrix $B$.
An alternative is to set $x  = B\inv y$ 
in \nref{eq:filtBase} in which case  the problem is restated as: 
\begin{equation} \label{eq:gen2}
K y = \rho(\lambda) B\inv y, \quad \mbox{with} \quad K=\rho(B\inv A)B\inv, \; y = B x 
\end{equation}
where again $K$ is symmetric and $B\inv $ is SPD.

Thus, we can apply Lanczos algorithms, along with a filter $\rho$ that is
constructed for  a given  interval, to  the matrix  pencil $(K,  B)$ in
\eqref{eq:gen1}  or  the pencil  $(K,  B\inv)$  in \eqref{eq:gen2}  to
extract  the desired  eigenvalues of  $(A,B)$ in  the interval.   When
applying  Algorithm~\ref{alg:LanAB}   to  \eqref{eq:gen1},   the  Ritz
vectors that  approximate eigenvectors of  $(A,B)$ take the  form $W_m
y_i$,  where $W_m  =  [w_1,\ldots, w_m]$  and  $(\theta_i,y_i)$ is 
eigenpair of  $T_m$.  On  the other hand,  for \eqref{eq:gen2}, 
 approximate  eigenvectors are of the form:
$B\inv W_m  y_i =  Z_m y_i$, where
$Z_m = [z_1,\ldots,  z_m]$.  Moreover, for the  polynomial filter, the
matrix $A$ should be first shifted and scaled as in \eqref{eq:mapping}
to ensure that the spectrum of $B\inv A$ is contained in $[-1,1]$.  In
terms  of computational  cost, the  difference between  performing the
Lanczos algorithm with \eqref{eq:gen1}  and with \eqref{eq:gen2} is small,
whereas for the rational filter, it might be often preferable to solve
\eqref{eq:gen2} with Algorithm \ref{alg:LanAB},  since solving
systems with $B$ can  be  avoided.  
A few practical  details  in applying  the  filters
follow.   First,   consider  applying  Algorithm   \ref{alg:LanAB}  to
\eqref{eq:gen1} with a polynomial filter  $\rho$.  
Line~3 of this algorithm 
can be stated  as follows
\begin{equation}
\label{eq:polyapp1}
 z:= B\rho(B\inv \hat A) w_i - \beta_i z_{i-1}
 = \rho(\hat A B\inv ) z_i - \beta_i z_{i-1},
\end{equation}
where $\hat A= (A-cB)/d$ as in \eqref{eq:mapping} and note that $z_i=Bw_i$.
So, applying the filtered matrix to a vector requires $k$ solves with $B$ and $k$ 
matrix-vector products with $A$, where $k$ is the degree of the polynomial filter.
Second, when rational filters are used, 
employing Algorithm~\ref{alg:LanAB}
 to the pair $(K, B\inv)$ in \eqref{eq:gen2}, 
Line 3 of this algorithm becomes
\begin{equation}
\label{eq:rationalapp}
 z := \rho(B\inv A)B\inv w_i - \beta_i z_{i-1} = \rho(B\inv A)z_i - \beta_i z_{i-1}, 
\end{equation}
where we have set    $w_i= Bz_i$.
Let the rational filter take the form  $\rho(t)= 2 \Re e \sum_j \sum_k a_{jk}(t-\sigma_j)^{-k}$.
The operation of $\rho(B\inv A)$ on a vector $z$ can be rewritten as
\begin{align} \label{eq:ratfiltapply}
\rho(B\inv A) z  &= 2 \Re e \sum_{j=1}^{p} \sum_{k=1}^{k_j} \alpha_{jk}  (B\inv A - \sigma_j I)^{-k} z
= 2 \Re e \sum_{j=1}^{p} \sum_{k=1}^{k_j} \alpha_{jk}  [(A-\sigma_j B)\inv B ]^{k} z ,
\notag 
\end{align}
which requires $\sum k_j$ matrix-vector products with $B$ 
and $\sum k_j$ solves with $(A-\sigma_j B)$,
where $k_j$ denotes  the number of the repetitions of the $j$-th pole.
The main operations involved in the application of the polynomial and 
rational filters for standard and generalized eigenvalue problems 
are summarized in Table~\ref{tab:itercomp}.

Third,  when  applying  Lanczos  algorithms  to  \eqref{eq:gen2},  the
initial vectors  $z_1$ and $w_1$  can be generated by  first selecting
$z_1$  as  a  random  vector  with  unit  $B$-norm  and  then  compute
$w_1=Bz_1$ in  order to  avoid a system solution  with $B$.  
 Finally, Lanczos
algorithms for  generalized eigenvalue problems  with the same  TR and
locking schemes as in Algorithm \ref{alg:FiltLanTR} have also been developed
in EVSL, where the deflation takes the form
\begin{equation}
(I-BUU\trans)Ax = \lambda Bx, \notag
\end{equation}
where $U$ has as its columns 
 the converged Ritz vectors, which are $B$-orthonormal.

\begin{table}[h]
\caption{The computations involving $A$ and $B$ in applying the filters, 
where $k$ denotes the degree of the polynomial filter and $k_j$ denotes 
the number of the repetition
of pole $\sigma_j$ of the
rational filter.
\label{tab:itercomp}}
\renewcommand{\arraystretch}{1.1}
\centering
\begin{tabular}{r|cc|cc}
 & \multicolumn{2}{c|}{standard} & \multicolumn{2}{c}{generalized}\tabularnewline
 & polynomial & rational & polynomial & rational\tabularnewline
\hline 
$Ax$ & $k$ &  & $k$ & \tabularnewline
$Bx$ &  &  &  & $\sum k_j$\tabularnewline
$(A-\sigma_j B)\inv x$ &  & $\sum k_j$ &  & $\sum k_j$ \tabularnewline
$B\inv x$ &  &  & $k$ & \tabularnewline
\end{tabular} 
\end{table}

Before closing this section, we point out that in addition to 
the Lanczos algorithm, EVSL also includes subspace iteration as
another projection method. Subspace iteration is generally slower
than Krylov subspace-based methods but it may have some advantages
in nonlinear eigenvalue problems where a subspace computed at a given 
nonlinear step can be used as initial subspace for  the next step.
Note that  subspace iteration  requires a  reasonable estimate  of the
number of the eigenvalues to  compute beforehand in order to determine
the subspace dimension. This number is readily available from the DOS
algorithm for spectrum slicing.

\subsection{Parallel EVSL}
\begin{figure}
\caption{The two main levels of parallelism in EVSL. \label{fig:parevsl}}

\vspace{0.5em}

\begin{center}
{\includegraphics[width=0.7\textwidth]{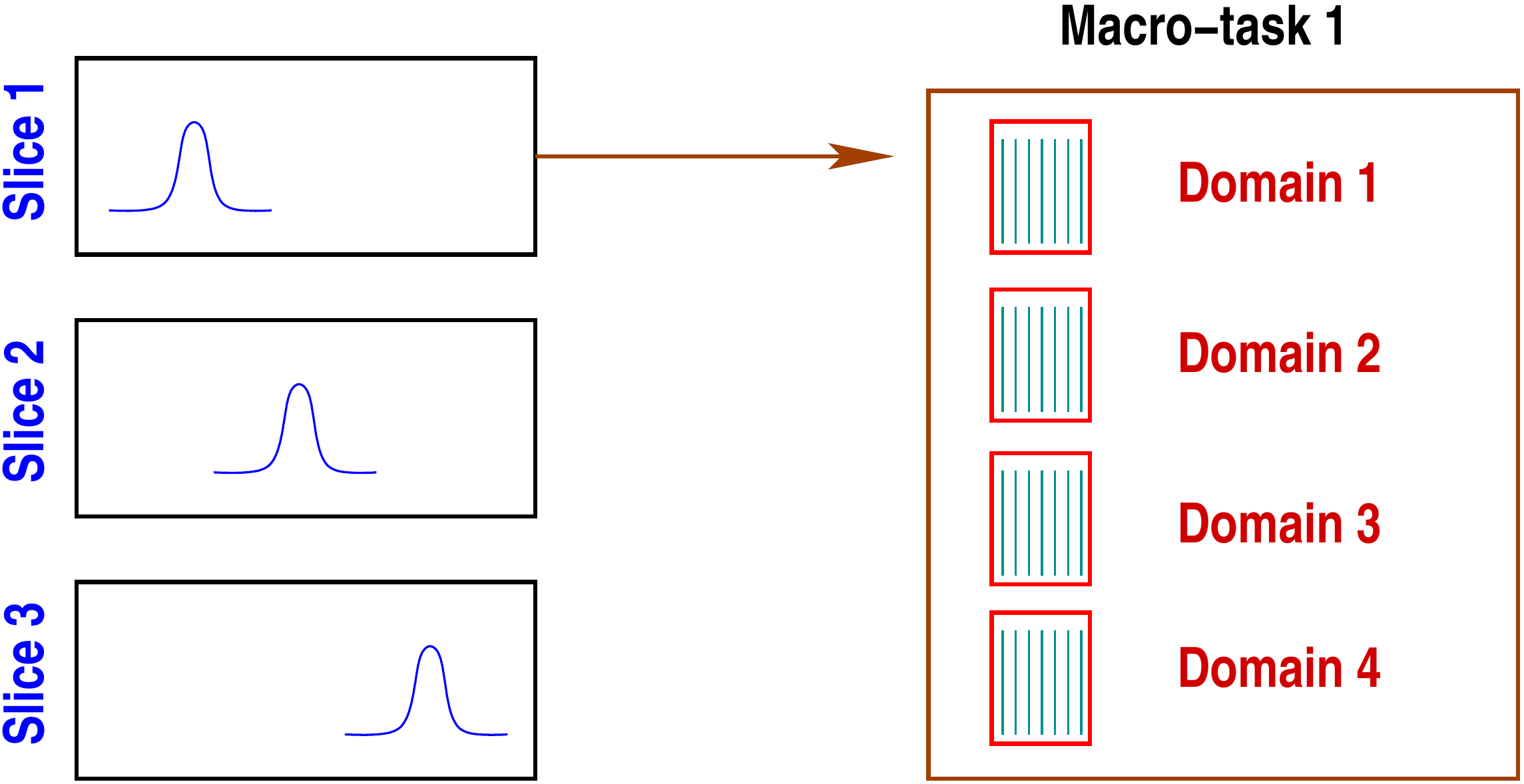}}
\end{center}
\end{figure}

This  paper  focuses  on  the  spectrum  slicing  algorithms  and  the
projection  methods with  filtering techniques  that are  included in
EVSL.  These algorithms  were implemented for  scalar machines
and  can  take   advantage  of standard multi-core computers with  shared  memory.
However, the algorithms in EVSL were designed with solving large-scale
problems on parallel platforms in  mind.  It should be clear that 
the  spectrum  slicing  framework   is  very  appealing  for  parallel
processing, where two main levels of parallelism can be exploited.  At
a macro level, we can divide the interval of interest into a number of
slices and map each slice to a group of processes, so that eigenvalues
contained  in different  slices can  be extracted  in parallel.   This
strategy forms the  first level of parallelism across  the slices.  On
parallel machines  with distributed memory, this  level of parallelism
can be easily implemented by organizing subsets of processes by groups
via subcommunicators in  a message passing model, such  as MPI.  
A second  level of parallelism exploits domain
decomposition across the processes within  a group, where each process
holds the data in its corresponding domain and the parallel matrix and
vector operations are performed  collectively within the group.  These
two levels  of parallelism  are sketched in  Figure \ref{fig:parevsl}.
More careful discussions on the parallel implementations, the parallel
efficiency and the  scalability of the proposed algorithms  are out of
the scope of this paper and will  be the object of a future article on
the  parallel   EVSL  package   (pEVSL),  which  is   currently  under
development.

\section{The EVSL package}
EVSL is a package that provides a collection of spectrum slicing algorithms
and projection methods for solving eigenvalue problems.
In this section, we address the software design and features of EVSL, 
its basic functionalities, and its user interfaces.

\subsection{Software design and user interface}
EVSL is written in C with an added Fortran interface and it is designed
primarily for large and sparse eigenvalue problems
in both the standard and the generalized forms, 
where the sparse matrices are stored in the CSR format.
The main functions of EVSL are listed next.
Two algorithms for computing spectral densities, namely, 
the kernel polynomial method \texttt{kpmdos} and the Lanczos algorithm \texttt{landos},
and the corresponding spectrum slicers, \texttt{spslicer} and \texttt{spslicer2} 
are included in EVSL.
Two filter options are available in EVSL, one based on Chebyshev polynomial expansions and the other one based on least-squares (L-S) 
rational approximations. These are built by the functions
\texttt{find\_pol} and \texttt{find\_ratf} respectively.
EVSL also includes 5 projection methods for eigenvalue problems:
\begin{itemize} 
\item Polynomial filtered non-restart Lanczos algorithm \texttt{ChebLanNr},
\item Rational filtered non-restart Lanczos algorithm \texttt{RatLanNr},
\item Polynomial filtered thick-restart Lanczos algorithm \texttt{ChebLanTr},
\item Rational filtered thick-restart Lanczos algorithm \texttt{RatLanTr}, and 
\item Polynomial filtered subspace iterations \texttt{ChebSI}.
\end{itemize} 
Two auxiliary functions named \texttt{LanBounds} and \texttt{LanTrbounds} 
implement simple Lanczos iterations and a more reliable thick-restart version
to provide sufficiently good estimations of the smallest and the largest eigenvalues,
which are needed to transform the spectrum onto $[-1,1]$ 
by the polynomial filtering approach.
A  subroutine  named   \texttt{lsPol}  computes  Chebyshev  polynomial
expansions that approximate functions $f(t)=t\inv$ and $f(t)=t^{-1/2}$
in a least-squares sense. This is used to approximate  the operations
$B\inv v$  and $B^{-1/2}v$   
as required  by the DOS algorithms  and the
polynomial filtering algorithm for generalized problems.  The accuracy
of  this  approximation   relies  on  the  condition   number  of  $B$
\cite{gdos}.     For   many    applications   from    finite   element
discretization, the matrix $B$ is the mass  matrix, which is
often very well-conditioned  after a proper diagonal scaling.  To
solve linear systems  with $A-\sigma_j B$, EVSL  currently relies on
sparse direct methods.  Along with the library that contains the above
mentioned algorithms, test programs are also provided to demonstrate the
use of  EVSL with  finite difference Laplacians  on regular  grids and
general matrices stored under the MatrixMarket format.

The algorithms in EVSL can also  be used in a `matrix-free' form that relies on
user-defined matrix-vector multiplication  functions and functions for
solving  the  involved  linear   systems.  These  functions  have  the
following  generic   prototypes  to   allow  common   interfaces  with
user-specific functions:
\begin{center}
\texttt{typedef void (*NAME)(double *x, double *y, void *data)} 
\end{center}
and
\begin{center}
\texttt{typedef void (*NAME)(double *x\_Re, double *x\_Im, double *y\_Re, double *y\_Im, void *data)} 
\end{center}
for the computations in real and complex arithmetics, where vector 
\texttt{x} is the input, vector \texttt{y} is the output and 
\texttt{data} points to all the required data to perform the function.
Below are two sample code snippets to illustrate the common user interfaces
of EVSL. The first example shows a call to
 \texttt{lanDos} to slice the spectrum followed by calls to
 the polynomial filtering with \texttt{ChebLanNr} to compute the 
eigenvalues of $(A,B)$ contained in an interval $(\xi,\eta)$ with CSR matrices
\texttt{Acsr} and \texttt{Bcsr}.
A sparse direct solver is first setup for 
the solve with $B$ and $L\trans$ where $B=LL\trans$ (lines 9-11).
Then, a lower and an upper bounds of the spectrum are obtained 
by \texttt{LanTrbounds} (line 11).
The DOS is computed by \texttt{lanDos} followed by the spectrum slicing
with \texttt{spslicer2}. Finally, we construct a polynomial filter for each
slice using \texttt{find\_pol} and combine the filter with
 \texttt{ChebLanNr} to extract the eigenvalues contained in the 
slice (lines 15-19).

\noindent
\begin{minipage}{\textwidth}
\begin{lstlisting}
/* Example 1: ChebLanNr for (A,B) */
#include "evsl.h"
EVSLStart();
SetGenEig();
SetAMatrix(&Acsr);
SetBMatrix(&Bcsr);
SetupBSolDirect(&Bcsr, &Bsol);
SetBSol(BSolDirect, Bsol);
SetLTSol(LTSolDirect, Bsol);
LanTrbounds(..., &lmin, &lmax);
intv[0] = xi; intv[1] = eta; intv[2] = lmin; intv[3] = lmax;
lanDos(..., npts, xdos, ydos, &nev_est, intv);
spslicer2(xdos, ydos, nslices, npts, sli);
for (i = 0; i < nslices; i++) {
  intv[0] = sli[i];  intv[1] = sli[i+1];
  find_pol(intv, &pol);
  ChebLanNr(intv, ..., &pol, &nev, &lam, &Y, &res, ...);
}
EVSLFinish();
\end{lstlisting}
\end{minipage}

The second example shows a code that  uses  rational filtering 
with a call to \texttt{RatLanNr}. The rational filter is built by 
\texttt{find\_ratf}. 
For each pole of the rational filter, we setup a  direct
solver for the corresponding shifted matrix and save the data
for all the solvers in \texttt{rat}.
Once the rational filter is properly setup, \texttt{RatLanNr} can be 
invoked with it.

\vspace{0.5em}

\noindent
\begin{minipage}{\textwidth}
\begin{lstlisting}
/* Example 2: RatLanNr for (A,B) */
for (i = 0; i < nslices; i++) {
  intv[0] = sli[i]; intv[1] = sli[i+1]; 
  intv[2] = lmin;   intv[3] = lmax;
  find_ratf(intv, &rat);
  void **sol = (void **)malloc(num_pole * sizeof(void *));
  SetupASIGMABSolDirect(&Acsr, &Bcsr, num_pole, rat, sol);
  SetASigmaBSol(&rat, ASIGMABSolDirect, sol);
  RatLanNr(intv, ..., &rat, &nev, &lam, &Y, &res, ...);
}
\end{lstlisting}
\end{minipage}

\subsection{Availability and dependencies}
EVSL is available from its website \url{http://www.users.cs.umn.edu/~saad/software/EVSL/index.html} and also
from the development website \url{https://github.com/eigs/EVSL}.
EVSL requires the dense linear algebra package BLAS \cite{Lawson:1979:BLA:355841.355847,Dongarra:1988:ESF:42288.42291,Dongarra:1990:SLB:77626.79170} and LAPACK \cite{laug}.
Optimized high performance BLAS and LAPACK routines, such as those
from Intel Math Kernel Library (MKL) and OpenBLAS \cite{Wang:2013:AAG:2503210.2503219},
and sparse matrix-vector multiplication routines
are  suggested for obtaining good performance.
For the linear system solves,
CXsparse \cite{doi:10.1137/1.9780898718881} is distributed along with EVSL
as the default sparse direct solver.
However, it is possible,  and we highly recommended to configure
EVSL  with other external high performance direct solvers, such as
Cholmod \cite{Chen:2008:ACS:1391989.1391995} and UMFpack 
\cite{Davis:2004:AUV:992200.992206} from SuiteSparse, and Pardiso
\cite{SCHENK200169} with the unified interface prototype.
Wrappers with this interface can be easily written for other solver 
options.

\section{Numerical experiments}\label{sec:num}
In this  section, we provide  a few examples to 
illustrate  the performance and ability of EVSL. In the tests we
  compute   interior  eigenvalues   of  Laplacian   matrices  and
Hamiltonian  matrices from electronic structure,  and we also solve
 generalized
eigenvalue  problems   from a  Maxwell  electromagnetics  problems.   The
experiments were  carried out on  one node  of Catalyst, a  cluster at
Lawrence  Livermore  National Laboratory,  which  is  equipped with  a
dual-socket Intel Xeon E5-2695 processor (24 cores) and 125 GB memory.
The EVSL package was compiled with the Intel \texttt{icc} compiler and
linked  with the  threaded Intel  MKL for  the
BLAS/LAPACK and the sparse matrix-vector multiplication routines.

\subsection{Laplacian matrices}
We begin  with discretized  Laplacians obtained  from the  5-point and
7-point finite-difference  schemes on regular  2-D and 3-D  grids.  In
each column  of Table~\ref{tab:lap},  we list  the grid sizes 
for each test case, the size of the matrix ($n$), the number of the nonzeros
($nnz$),  the  spectral interval ($[a,b]$), two  intervals  of
interest  from which  to extract  eigenvalues ($[\xi,\eta]$),  and the
actual numbers  of the eigenvalues contained  in $[\xi,\eta]$, denoted
by $\nu_{[\xi,\eta]}$.
\begin{table}[tbh]
\caption{Discretized 5-pt/7-pt Laplacians  on 2-D/3-D regular grids. \label{tab:lap}}%
\centering
\def\arraystretch{1.1}
\begin{tabular}
[c]{r|rrrcc}
\multicolumn{1}{c|}{Grid} & \multicolumn{1}{c}{n} & \multicolumn{1}{c}{nnz} & \multicolumn{1}{c}{$[a,b]$} & $[\xi,\eta]$& $\nu_{[\xi,\eta]}$ \\\hline
\multirow{2}{*}{$343^2$} & \multirow{2}{*}{$117,649$} & \multirow{2}{*}{$586,873$} & \multirow{2}{*}{$[0,7.9998]$} & $[0.40, 0.436]$ & 356 \\
& & & & $[1.00, 1.033]$ & 347\\
\hline
\multirow{2}{*}{$729^2$} & \multirow{2}{*}{$531,441$} & \multirow{2}{*}{$2,654,289$} & \multirow{2}{*}{$[0,7.9998]$} & $[0.40,0.410]$& 446 \\
& & & & $[1.00,1.009]$& 435\\
\hline
\multirow{2}{*}{$1000^2$} & \multirow{2}{*}{$1,000,000$} & \multirow{2}{*}{$4,996,000$} & \multirow{2}{*}{$[0,7.9998]$} & $[0.40,0.405]$& 429 \\
& & & & $[1.00,1.005]$& 460\\
\hline
\hline
\multirow{2}{*}{$49^3$} & \multirow{2}{*}{$117,649$} & \multirow{2}{*}{$809,137$} & \multirow{2}{*}{$[0,11.9882]$} & $[0.40, 0.570]$ & 343\\
& & & & $[1.00,1.100]$& 345 \\
\hline
\multirow{2}{*}{$81^3$} & \multirow{2}{*}{$531,441$} & \multirow{2}{*}{$3,680,721$} & \multirow{2}{*}{$[0,11.9882]$} & $[0.40, 0.450]$ & 433\\
& & & & $[1.00,1.028]$& 420\\
\hline
\multirow{2}{*}{$100^3$} & \multirow{2}{*}{$1,000,000$} & \multirow{2}{*}{$6,940,000$} & \multirow{2}{*}{$[0,11.9882]$} & $[0.40, 0.428]$ & 454\\
& & & & $[1.00,1.018]$& 475\\
\end{tabular}
\end{table}

Table~\ref{tab:lap_pol_results} shows test results for computing the
eigenvalues  of  the  Laplacians  in   the  given  intervals  and  the
corresponding  eigenvectors using  \texttt{ChebLanNr}, where  `deg' is
the degree  of the polynomial  filter, `iter' indicates the  number of
Lanczos  iterations   and  `nmv'   stands  for   the  number   of 
matrix-vector products required.  For the timings, we  report the time
for  computing  the  matrix-vector  products (`t-mv'),  the  time  for
performing the  reorthogonalization (`t-orth'), and the  total time to
solution (`t-tot').  As shown, the  interval that is deeper inside the
spectrum requires a higher-degree  polynomial filter and
more iterations  as well to  extract all the eigenvalues, in  general.   This makes
the computation  more expensive than  with the other  interval for
computing  roughly the  same number  of eigenvalues.   This issue  was
known from the results in  \cite{spectrumslicing}. For the two largest
2-D and  3-D problems,  the total  time was dominated  by the  time of
performing the matrix-vector products.

\begin{table}[tbh]
\caption{Polynomial filtered Lanczos algorithm for discretized Laplacians. 
Times are measured in seconds.}
\label{tab:lap_pol_results}%
\def\arraystretch{1.1}
\begin{tabular}[c]{r|crcrrrr}
\multicolumn{1}{c|}{\multirow{2}{*}{Grid}}  &
\multicolumn{1}{c}{\multirow{2}{*}{$[\xi,\eta]$} }  &
 \multicolumn{1}{c}{\multirow{2}{*}{deg}} & \multirow{2}{*}{niter} & 
\multicolumn{1}{c}{\multirow{2}{*}{nmv}}
 & \multicolumn{3}{c}{iter}  \\  
   & &  & &  & \multicolumn{1}{c}{t-mv} & t-orth & \multicolumn{1}{c}{t-tot}  
\\\hline    
\multirow{2}{*}{$343^2$}&$[0.40,0.4360]$ & 157 & 1,290 & 203,200 & 10.66 & 24.11 & 51.95 \\
&$[1.00, 1.0330]$ & 256 &  1,310 & 336,219 & 17.39 & 25.48 & 66.71  \\ \hline
\multirow{2}{*}{$729^2$}&$[0.40, 0.4100]$ &  557 &  1,610 & 898,330 & 182.51 & 168.20 & 497.78 \\
&$[1.00, 1.0090]$ & 936 & 1,590 & 1,490,547 & 294.33 & 179.45 & 679.08 \\\hline
\multirow{2}{*}{$1000^2$}&$[0.40, 0.4050]$ & 1,111  & 1,530   &  1,702,481 &  1451.40 & 308.98 & 2410.99  \\
&$[1.00, 1.0050]$ & 1,684 & 1,670  & 2,816,108 & 2402.99 & 367.49 & 3804.91  \\\hline\hline
\multirow{2}{*}{$49^3$}&$[0.40,0.5700]$ & 43 & 1,290 & 55,899  & 3.81  & 26.53  & 41.03  \\ 
&$[1.00,1.1000]$ & 107 & 1,270  & 136,449  & 9.80 & 24.25 & 47.00  \\\hline
\multirow{2}{*}{$81^3$}&$[0.40, 0.4500]$ & 141 & 1,590 & 224,905 & 77.24 & 189.91 & 359.61 \\
&$[1.00,1.0280]$ &  378 & 1,570 &  594,636   & 191.42  & 168.33 & 494.66  \\ \hline
\multirow{2}{*}{$100^3$}&$[0.40, 0.4280]$ & 248 &  1,710 & 425,030 & 563.36 & 398.70 &  1250.26  \\
&$[1.00,1.0180]$ & 588  & 1,830  & 1,077,689  & 1357.12  & 426.17  & 2314.43  \\ 
\end{tabular}
\end{table}

Table~\ref{tab:lap_rat_results}  presents  test   results  when  using
\texttt{RatLanNr}   to    solve   the    same   problems    given   in
Table~\ref{tab:lap}.  We  used the  sparse direct solver  Pardiso from
the  Intel MKL  to solve  the  complex symmetric  linear systems  when
applying the rational filter.  Comparing the CPU times (t-tot) in this
table with  those in Table~\ref{tab:lap_pol_results}, we  can see that
it is  more efficient to  use the rational filtered  Lanczos algorithm
for the  2-D problems than  the polynomial filtered  counterpart.  The
numbers in  the columns  labeled `fact-fill'  and `fact-time'  are the
fill-factors (defined  as the ratio of  the number of nonzeros  of the
factors over  the number of nonzeros  of the original matrix)  and the
CPU times for  factoring the shifted matrices,  respectively.  For the
2-D problems,  the fill-factors  stay low  and the  factorizations are
inexpensive.   On the  other hand,  for the  3-D problems,  the sparse
direct solver  becomes much more  expensive: it  is not only  that the
factorizations become much more costly, as reflected by the higher and
rapidly  increasing fill-factors,  but  that the  solve  phase of  the
direct  solver,  performed  at   each  iteration,  becomes  more  time
consuming as well. The timings for  performing the solves are shown in
the  column  labeled  `t-sv'.   Consequently,  for the 3-D problems
the  rational  filtered Lanczos algorithm is less efficient than 
the  polynomial  filtered   algorithm.   Finally,  it  is  worth
mentioning that the rational filtered Lanczos algorithm required far fewer
iterations than  the polynomial filtered algorithm  (about half). This
indicates that the  quality of the rational filter is  better than the
polynomial  filter with  the default  threshold.  Rational  filters are
much more effective  in amplifying the eigenvalues in  a target region
while damping the unwanted eigenvalues.

\begin{table}[tbh]
\caption{Rational filtered Lanczos algorithm for discretized Laplacians.
Times are measured in seconds.
\label{tab:lap_rat_results}}%
\centering
\def\arraystretch{1.1}
\begin{tabular}
[c]{r|cccrrrrr}
\multicolumn{1}{c|}{\multirow{2}{*}{Grid}}  
&\multicolumn{1}{c}{\multirow{2}{*}{$[\xi,\eta]$}}  & \multirow{2}{*}{niter} & \multirow{2}{*}{nsv} & 
\multicolumn{2}{c}{fact} & \multicolumn{3}{c}{iter} \\
  &  & & & \multicolumn{1}{c}{fill} & \multicolumn{1}{c}{time} & \multicolumn{1}{c}{t-sv} & t-orth & 
\multicolumn{1}{c}{t-tot} \\\hline
\multirow{2}{*}{$343^2$}&$[0.40,0.4360]$ & 590 & 1,184 & \multirow{2}{*}{10.3} & \multirow{2}{*}{0.62} & 39.23  & 6.11 & 50.45   \\
&$[1.00, 1.0330]$ & 590 & 1,184 &  &  & 38.13  & 5.23 & 47.66  \\
\hline
\multirow{2}{*}{$729^2$}&$[0.40, 0.4100]$ & 730 & 1,464 & \multirow{2}{*}{13.0} & \multirow{2}{*}{2.50} & 274.70 & 37.86 & 339.02  \\
&$[1.00,1.0090]$& 710   & 1,424 &  &  & 269.14 & 37.38 & 331.58 \\
\hline
\multirow{2}{*}{$1000^2$}&$[0.40, 0.4050]$ & 710 & 1,424 & \multirow{2}{*}{14.1} & \multirow{2}{*}{4.77} & 534.90  & 75.46  & 659.40 \\
&$[1.00,1.0050]$& 750  & 1,504 &   &  & 559.76  &  86.11 & 706.61 \\
\hline\hline
\multirow{2}{*}{$49^3$} & $[0.40,0.5700]$ & 630  & 1,264 & \multirow{2}{*}{68.6} & \multirow{2}{*}{1.90}  & 168.56  & 6.43  & 180.11 \\
&$[1.00,1.1000]$ & 610  & 1,224 &   &  & 159.86  & 5.97  & 170.63 \\
\hline 
\multirow{2}{*}{$81^3$}&$[0.40,0.4500]$ & 750  & 1,504 & \multirow{2}{*}{149.7} & \multirow{2}{*}{28.77}  &   1850.33 & 42.14 & 1923.36\\
&$[1.00,1.0280]$  &  730  & 1,464 &  &   &   1813.15 & 40.40 & 1882.21 \\
\hline
\multirow{2}{*}{$100^3$}&$[0.40,0.4280]$& 810  &  1,624 & \multirow{2}{*}{186.3}&  \multirow{2}{*}{87.64} &  4343.00 & 104.30 & 4510.28 \\
&$[1.00,1.0180]$  & 930 & 1,864 &  &  & 5035.75 & 127.28 & 5240.53 \\
\end{tabular}
\end{table}

\subsection{Spectral slicing}

In  this  set of  experiments,  we  show  performance results  of  the
spectrum slicing algorithm with KPM in EVSL, where we computed all the
$1,971$ eigenvalues in  the interval $[0,1]$ of  the 7-point Laplacian
matrix on the 3-D grid of size  $49 \times 49 \times 49$. The interval
$[0,1]$ was first partitioned into up to $6$ slices in such a way that
each  slice contains  roughly the  same number  of eigenvalues.  Then,
\texttt{ChebLanNr} was used  to extract the eigenvalues  in each slice
individually.     The    column    labeled    $\nu_{[\xi,\eta]}$    in
Table~\ref{tab:lapslicing} lists  the number of  eigenvalues contained
in  the divided  subintervals, which  are fairly  close to  each other
across the slices.  For the  iteration time `t-tol', since on parallel
machines,  a parallelized  EVSL code  can compute  the eigenvalues  in
different slices  independently, the total  time of computing  all the
eigenvalues in the  whole interval will be the  maximum iteration time
across all the slices.  As shown,  compared with the CPU time required
by the solver with a single  slice, significant CPU time reduction can
be achieved with multiple slices.

With regard to load balancing,  the memory requirements were very well
balanced across the slices, since the memory allocation in the Lanczos
algorithm is  proportional to  the number  of eigenvalues  to compute.
Note also that the iteration times shown in the column labeled `t-tol'
were also reasonably  well balanced, except for the  first slice which
is a  ``boundary slice'' on the  left end of the  spectrum.  A special
type  of  polynomial filters  were  used  for boundary  slices,  which
usually  have lower  degrees than  with the  filters for  the internal
slices.  Moreover, the convergence rates for computing the eigenvalues
for boundary slices are typically  better, as reflected by the smaller
number of  iterations.  Therefore,  the total  iteration time  for the
first slice is much smaller than  with the other slices. Starting with
the  second  slice,   the  overall  iteration  time   did  not  change
dramatically.  This is  in spite of the fact that  the required degree
of the  polynomial filter keeps  increasing as the slice  moves deeper
inside the spectrum leading to a  higher cost for applying the filters
each time.   The explanation  is that  in these  problems the  cost of
performing matrix-vector  products was  insignificant relative  to the
cost of the reorthogonalizations.

\begin{table}[tbh]
\caption{Polynomial filtered  Lanczos algorithm for computing $1,971$ eigenvalues of a Laplacian of size $49^3$ within $[0, 1]$. \label{tab:lapslicing}}%
\def\arraystretch{1.1}
\begin{tabular}
[c]{c|crrrrrrr}
\multicolumn{1}{c|}{\multirow{2}{*}{Slices}}&\multicolumn{1}{c}{\multirow{2}{*}{$[\xi,\eta]$}}  &\multicolumn{1}{c}{\multirow{2}{*}{$\nu_{[\xi,\eta]}$}}  & \multirow{2}{*}{deg} & \multirow{2}{*}{niter} & \multicolumn{1}{c}{\multirow{2}{*}{nmv}} & \multicolumn{3}{c}{iter}  \\
& &  & &  & & t-mv & t-orth & \multicolumn{1}{c}{t-tot}  
\\\hline 
1&$[0.00000, 1.00000]$ &1,971 & 5 &  4,110 & 22,531 & 2.07 & 224.97 &  353.09  \\\hline
\multirow{2}{*}{2}
&$[0.00000,0.65863]$ & 997 & 6 & 2,230 & 14,389 & 1.34  & 69.92 & 107.76  \\ 
&$[0.65863, 1.00000]$ & 974 & 28 &  3,470 & 98,190 & 7.06 & 173.09 & 240.64  \\\hline
\multirow{3}{*}{3}
&$[0.00000,0.51208]$ & 657 & 7 & 1,510 &  11,241  & 1.06  & 34.43  & 52.09   \\ 
&$[0.51208,0.78294]$ & 662 & 31 & 2,390  & 74,814  & 5.66 &  88.42 & 126.46   \\ 
&$[0.78294,1.00000]$ & 652 & 46 & 2,350  & 108,844  & 7.79 &  76.49 & 116.00 \\\hline
\multirow{4}{*}{4}
&$[0.00000, 0.42983]$ & 495 & 8& 1,150 & 9,711 & 0.92 & 20.10 & 30.23   \\ 
&$[0.42983, 0.65394]$ & 484 & 35 & 1,770 & 62,504 & 4.71 & 44.65 & 68.66  \\ 
&$[0.65394, 0.83844]$ & 502 & 49 & 1,850 & 91,250 & 6.07 & 45.80 & 70.80  \\ 
&$[0.83844,1.00000]$ & 490 & 62 & 1,850  & 115,314 & 7.65 & 48.38 & 76.45 \\\hline
\multirow{5}{*}{5} 
&$[0.00000,0.37793]$& 386 & 8 &  970 & 8,162 & 0.73 & 15.34 & 22.46  \\ 
&$[0.37793,0.57373]$ & 401 & 37 & 1,450 & 54,125 & 4.00 & 29.74 & 46.00   \\ 
&$[0.57373, 0.73473]$ & 399 & 53 & 1,470  & 78,415 & 5.73 & 29.43 & 48.69   \\ 
&$[0.73473, 0.87424]$ & 400 & 68 & 1,450  & 99,136 & 7.36 & 29.89 & 51.41   \\ 
&$[0.87424,1.00000]$ & 385 & 81 & 1,430 & 116,377 & 8.83  & 29.51 & 53.57  \\\hline
\multirow{6}{*}{6} 
&$[0.00000,0.33926]$ & 329 & 9 &  810 & 7,637 & 0.68 & 13.30 & 20.84  \\ 
&$[0.33926,0.51429]$ & 328 & 40 & 1,212  &  48,808 & 3.68 & 23.85 & 37.20   \\ 
&$[0.51429, 0.65913]$ & 340 & 56 & 1,230  & 69,332 &  4.96 & 25.57 & 41.03  \\  
&$[0.65913,0.78384]$ & 322 & 72 & 1,230  & 89,026  & 6.29 & 25.65 & 44.98   \\ 
&$[0.78384,0.89719]$ & 345 & 85 & 1,230  & 105,065  & 7.42 & 26.21 & 46.52  \\ 
&$[0.89719,1.00000]$ & 307 & 100 & 1,190 & 119,507 & 8.13  & 23.93 & 43.53 \\
\end{tabular}
\end{table}

\subsection{Hamiltonian matrices}
In this  set of experiments,  we computed eigenpairs of  5 Hamiltonian
matrices from  the Kohn-Sham  equations for density  functional theory
calculations.  These matrices  are available from the  PARSEC group in
the             SuiteSparse              Matrix             Collection
\cite{Davis:2011:UFS:2049662.2049663}.   The  matrix   size  $n$,  the
number of  the nonzeros $nnz$ and  the range of the  spectrum $[a, b]$
are  provided in  Table  \ref{tab:Hamilton}.  Each  Hamiltonian has  a
number of occupied states of the  chemical compound and this number is
often  denoted  by $n_0$.  In  Density  Functional Theory  (DFT),  the
$n_0$-th  smallest  eigenvalue corresponds  to  the  Fermi level.   As
in~\cite{Filtlan-paper}, we computed all  the eigenvalues contained in
the interval  $[0.5n_0, 1.5n_0]$  and the  corresponding eigenvectors.
The target intervals, which are denoted by \intv{}, and the numbers of
the eigenvalues inside the intervals are given in the last two columns
of Table \ref{tab:Hamilton}.

\begin{table}[tbh]
\caption{
Hamiltonian matrices from the 
PARSEC set 
}%
\label{tab:Hamilton}
\def\arraystretch{1.1}
\begin{tabular}
[c]{l|crccc}
\multicolumn{1}{c|}{Matrix} & n & \multicolumn{1}{c}{nnz} & $[a,b]$ & $[\xi,\eta]$& $\nu_{[\xi,\eta]}$ \\\hline
$\mathrm{Ge_{87}H_{76}}$ & $112,985$ & $7,892,195$ & $[-1.214, 32.764]$ & $[-0.645, -0.0053]$ & 212 \\
$\mathrm{Ge_{99}H_{100}}$ & $112,985$ & $8,451,295$ & $[-1.226,32.703]$& $[-0.650, -0.0096]$&250 \\
$\mathrm{Si_{41}Ge_{41}H_{72}}$& $185,639$ & $15,011,265$ & $[-1.121,49.818]$&  $[-0.640,-0.0028]$& 218 \\
$\mathrm{Si_{87}H_{76}}$ &$240,369$ & $10,661,631$ & $[-1.196,43.074]$& $[-0.660,-0.0300]$& 213 \\
$\mathrm{Ga_{41}As_{41}H_{72}}$ & $268,096$ & $18,488,476$ & $[-1.250,1300.9]$&$[-0.640,-0.0000]$& 201\\
\end{tabular}
\end{table}


We  first report  in  Table~\ref{tab:Parsecresultpol} the  performance
results  of \texttt{ChebLanNr}  for computing  the eigenvalues  of the
Hamiltonians in  the target interval.   For the first 4  matrices, the
polynomial degrees used in the filters  are relatively low so that the 
cost of the matrix-vector products is not high relative to
  the  cost  of the  reorthogonalizations.   Extracting
eigenvalues  from the  last matrix  was more  challenging because  its
spectrum  is much  wider  than  the others,  with  the  result that  a
polynomial  filter of  a much  higher degree  was needed resulting in
a much higher cost for  matrix-vector products.  Similar results
were also reported in \cite{Filtlan-paper}.

\begin{table}[tbh]
\caption{Numerical results of polynomial filtered  Lanczos algorithm for PARSEC matrices. \label{tab:Parsecresultpol}}%
\centering
\def\arraystretch{1.1}
\begin{tabular}
[c]{l|rrrrrr}
\multicolumn{1}{c|}{\multirow{2}{*}{Matrix}}  & \multirow{2}{*}{deg} & \multirow{2}{*}{niter} & 
\multicolumn{1}{c}{\multirow{2}{*}{nmv}} & \multicolumn{3}{c}{iter}  \\
 &  & &  & \multicolumn{1}{c}{t-mv} & t-orth & \multicolumn{1}{c}{t-tot}  \\\hline
$\mathrm{Ge_{87}H_{76}}$ & 26 & 990 &  26,004 & 28.53 & 14.88 & 50.05  \\ 
$\mathrm{Ge_{99}H_{100}}$ & 26 & 1,090 & 28,642 & 33.71 & 17.52  & 59.18   \\ 
$\mathrm{Si_{41}Ge_{41}H_{72}}$ & 32 & 950  & 30,682  & 71.61 & 28.82 & 112.11  \\ 
$\mathrm{Si_{87}H_{76}}$ & 29 & 1,010 & 29,561 & 46.38  & 30.27 & 89.92  \\ 
$\mathrm{Ga_{41}As_{41}H_{72}}$ & 172 & 910  & 157,065  & 497.39 & 29.25 & 558.26 \\
\end{tabular}
\end{table}


Next, we present the performance results 
of \texttt{RatLanNr} in Table~\ref{tab:Parsecresultrat}.
An issue we confronted when using this algorithm for this set of matrices is
that the factorization of the shifted matrix becomes prohibitively expensive,
in terms of both the memory requirement as indicated by the 
high fill-ratios (e.g., the factorization of the last matrix required
42 GB memory), and the high factorization time
(remarkably, this  is even higher than the total iteration time 
with the polynomial filter).
Despite requiring only about half of the number of iterations that were needed
by \texttt{ChebLanNr}, the total iteration time using
the rational filter was still much higher, since the solve phase of 
the direct solver was also much more expensive compared with the matrix-vector products
performed in the polynomial filtered Lanczos algorithm.

\begin{table}[tbh]
\caption{Numerical results of rational filtered Lanczos algorithm for PARSEC matrices. \label{tab:Parsecresultrat}}%
\centering
\def\arraystretch{1.1}
\begin{tabular}
{l|ccrrrrr}
\multicolumn{1}{c|}{\multirow{2}{*}{Matrix}}  & \multirow{2}{*}{niter} & \multirow{2}{*}{nsv} & \multicolumn{2}{c}{fact} & \multicolumn{3}{c}{iter} \\
  &  & & \multicolumn{1}{c}{fill} & \multicolumn{1}{c}{time} & \multicolumn{1}{c}{t-sv} & t-orth & 
\multicolumn{1}{c}{t-tot} \\
\hline
  $\mathrm{Ge_{87}H_{76}}$ & 450 & 904 & 149.2 & 120.56 & 2022.33 & 4.03 & 2029.11  \\
$\mathrm{Ge_{99}H_{100}}$ & 490 & 984 & 146.5 & 152.07 & 2259.01 & 4.01 & 2266.42  \\
$\mathrm{Si_{41}Ge_{41}H_{72}}$ & 430 & 864  & 178.9  & 447.89 & 4587.51 & 7.09 & 4599.10 \\
$\mathrm{Si_{87}H_{76}}$ & 450 & 904 &  376.9  &  718.00  & 7058.93 & 9.71 & 7074.30 \\
$\mathrm{Ga_{41}As_{41}H_{72}}$ & 410 & 824  & 378.5  & 865.43 & 7456.54 & 9.97 & 7472.38 \\
\end{tabular}
\end{table}

\subsection{Generalized eigenvalue problems: the Maxwell Eigenproblem}
\begin{figure}[tbh]
\caption{Spectrum of the Maxwell eigenproblem on a 2-D square-disc and on a 3-D cube\label{fig:maxwellspectrum}}
\hspace{-2em}
\includegraphics[width=0.58\textwidth]{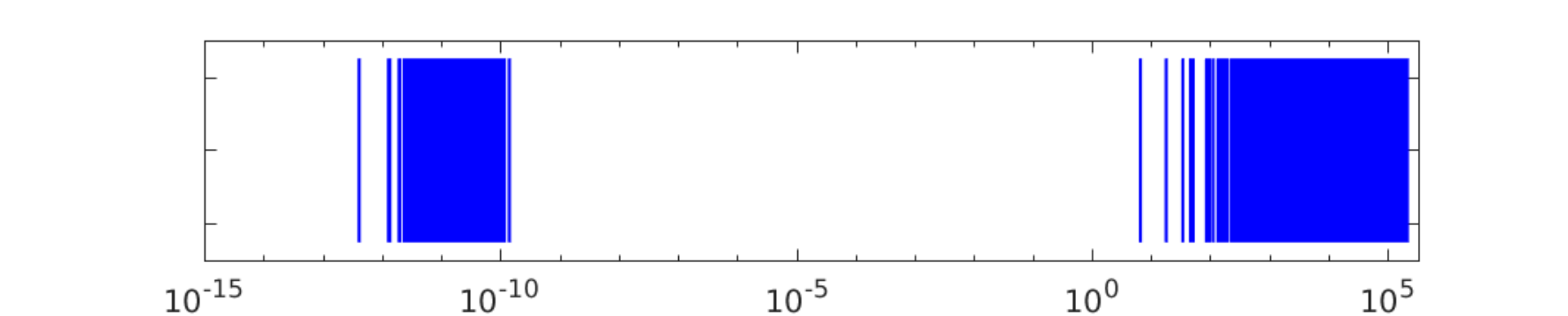}
\hspace{-3em}
\includegraphics[width=0.58\textwidth]{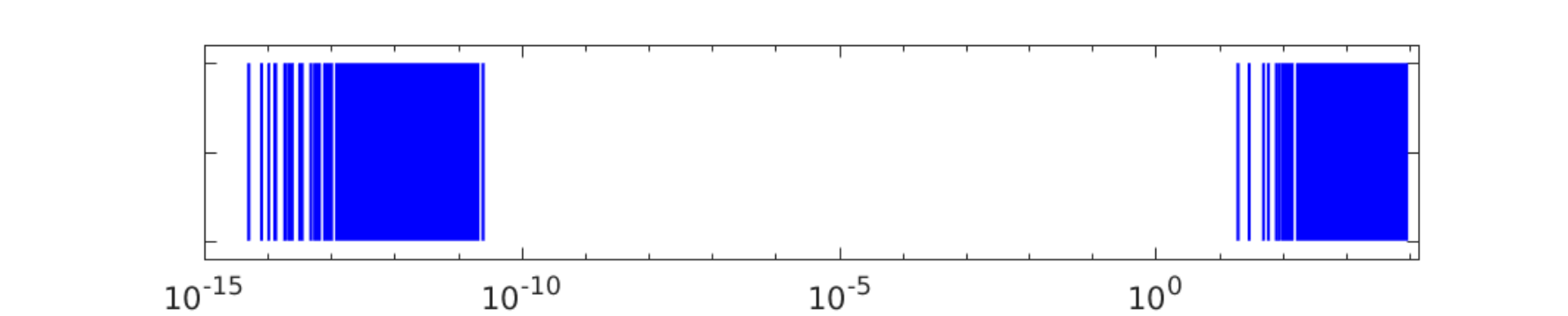} 
\end{figure}

\begin{figure}[tbh]
\caption{An illustration of the 1st, 10th, 18th, 29th, 49th, 59th, 65th and 75th electromagnetic eigenmodes of the Maxwell eigenproblem on 
a 2-D square-disc and on a 3-D cube \label{fig:eigmodes}}

\vspace{0.5em}

\includegraphics[width=0.1\textwidth]{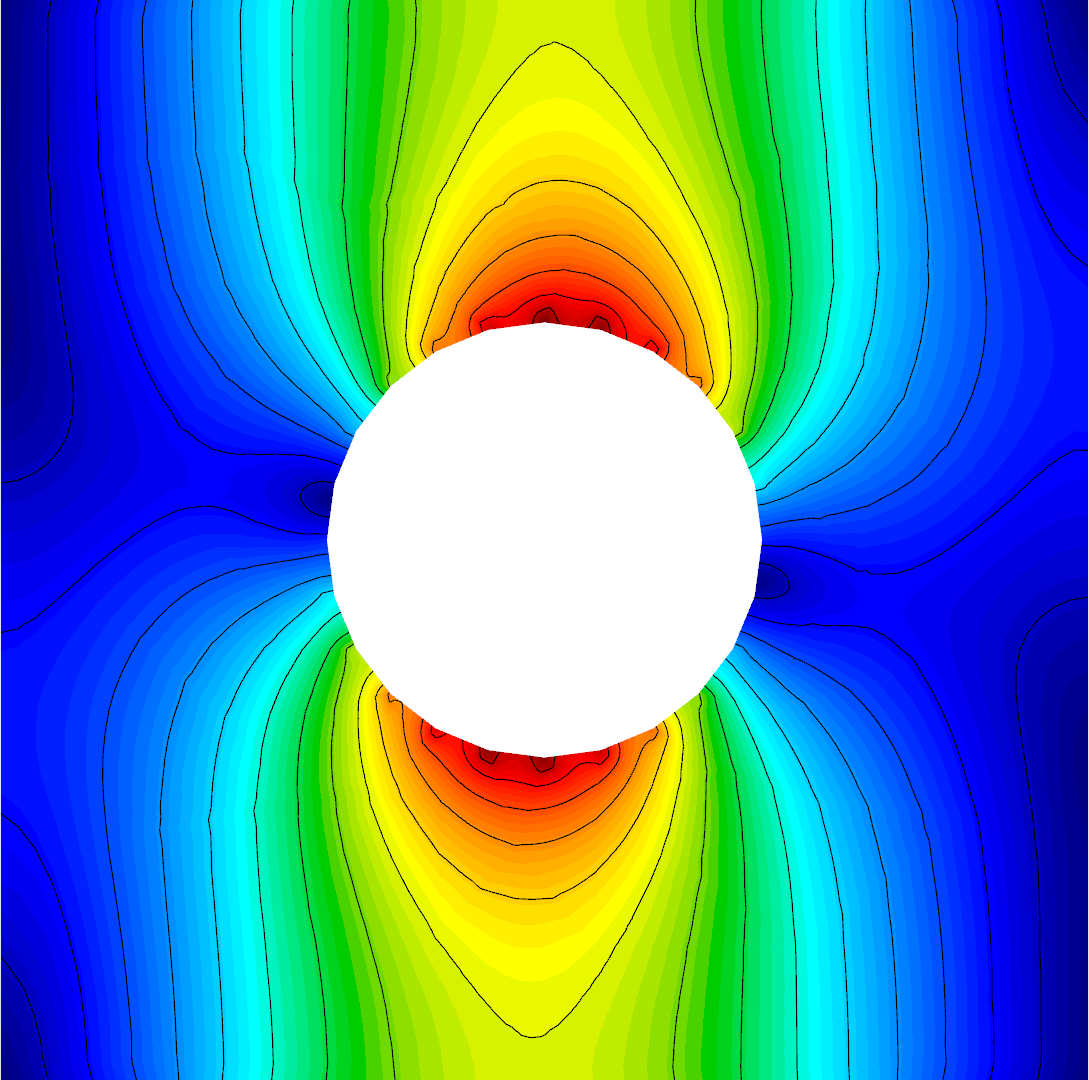} \hspace{0.2em}
\includegraphics[width=0.1\textwidth]{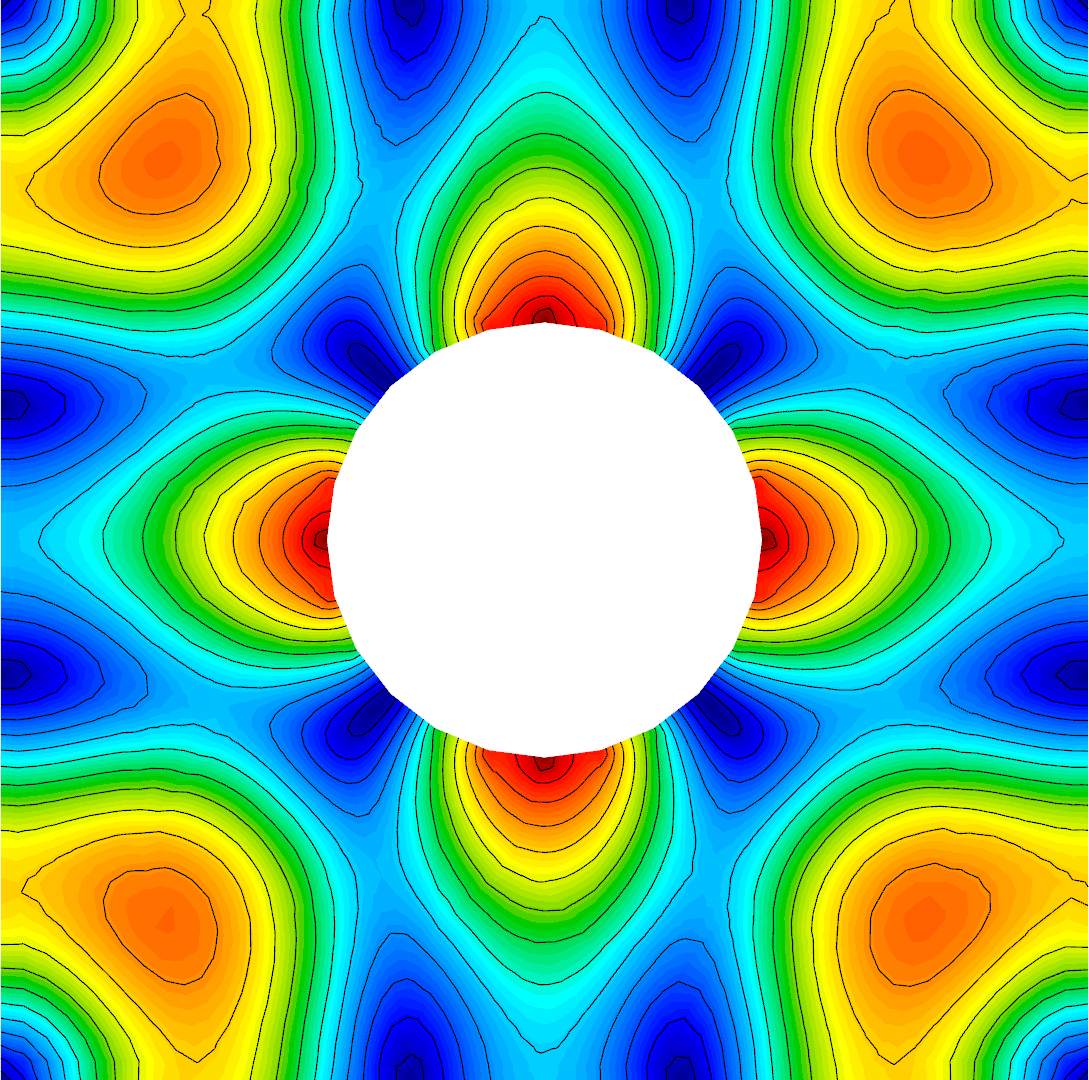} \hspace{0.2em}
\includegraphics[width=0.1\textwidth]{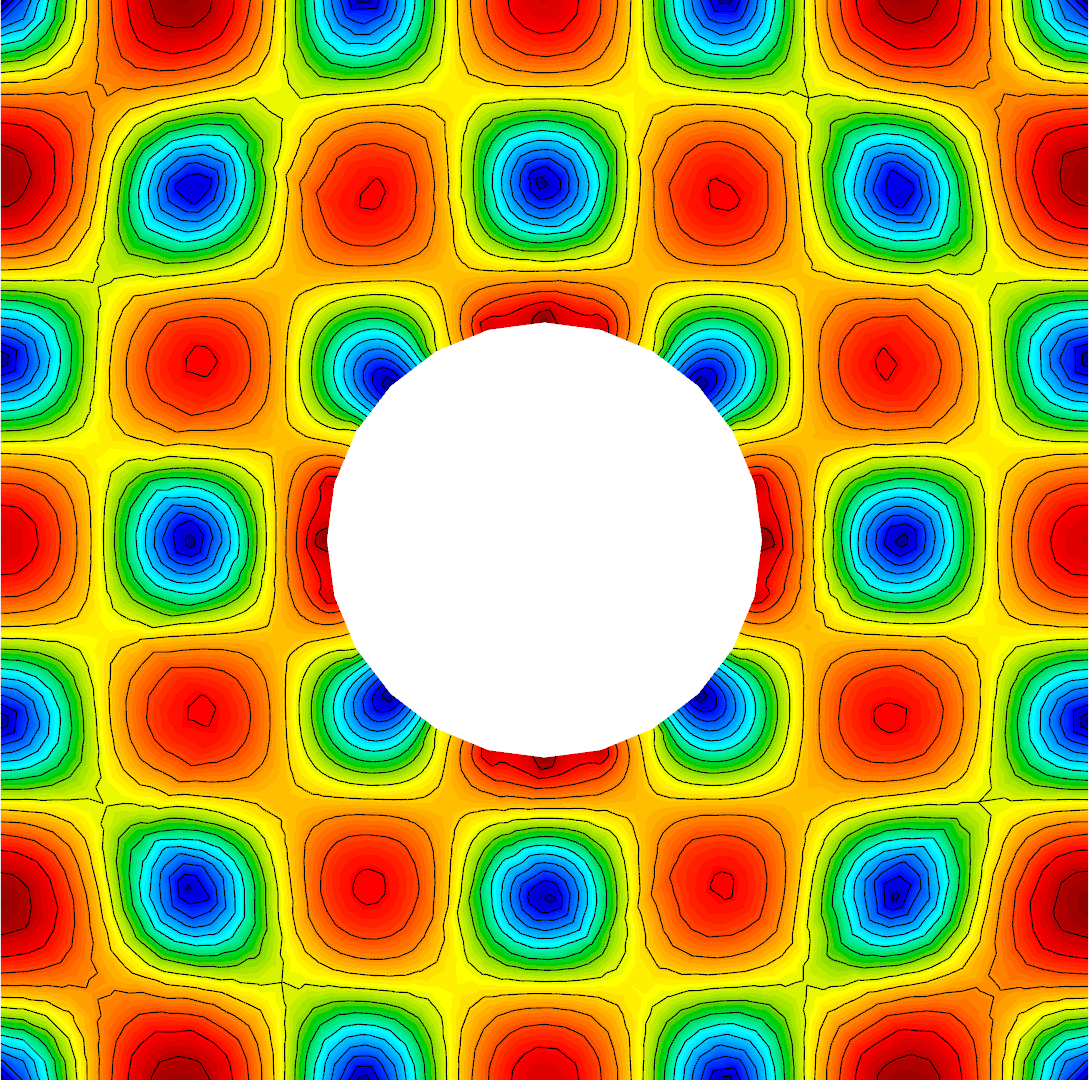} \hspace{0.2em}
\includegraphics[width=0.1\textwidth]{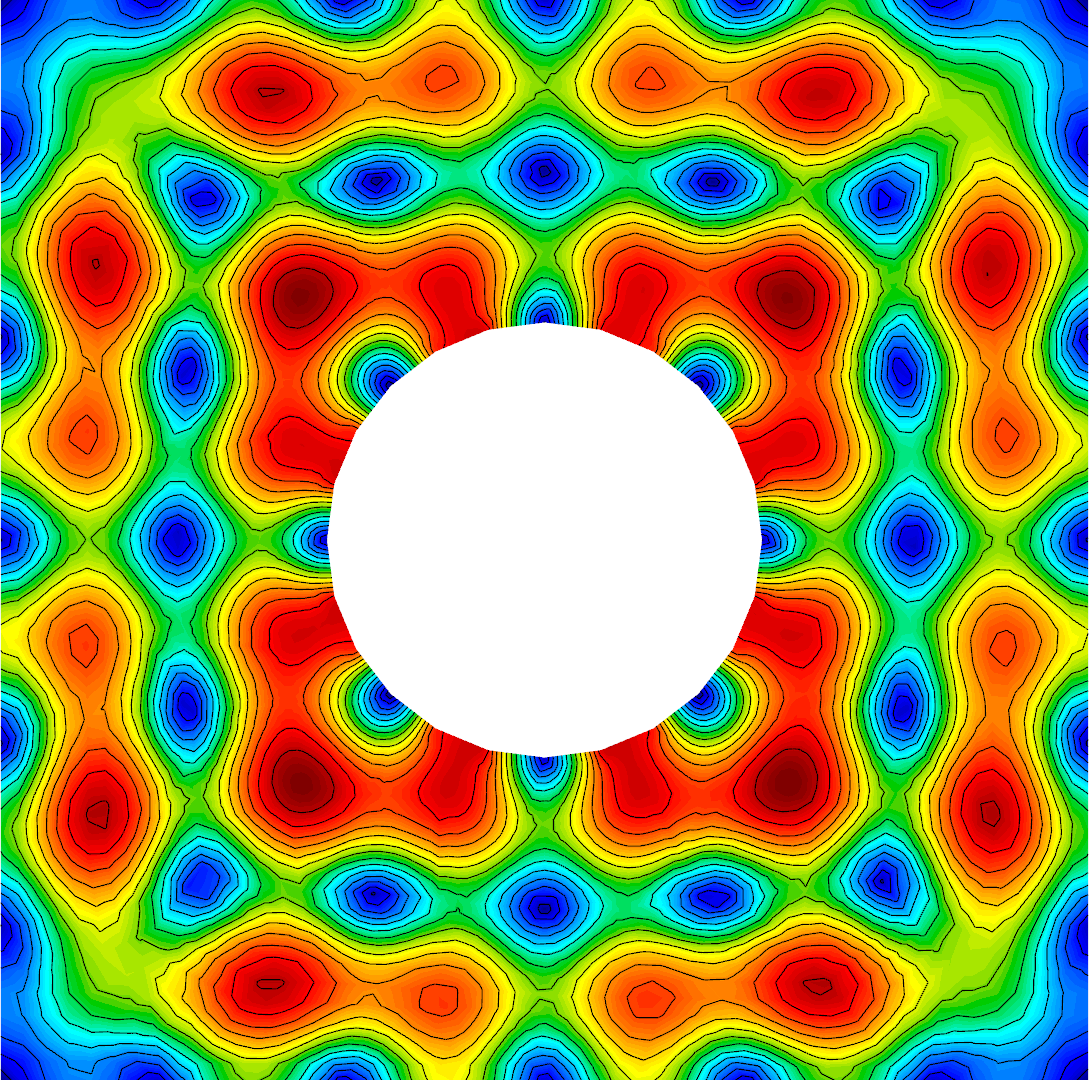} \hspace{0.2em}
\includegraphics[width=0.1\textwidth]{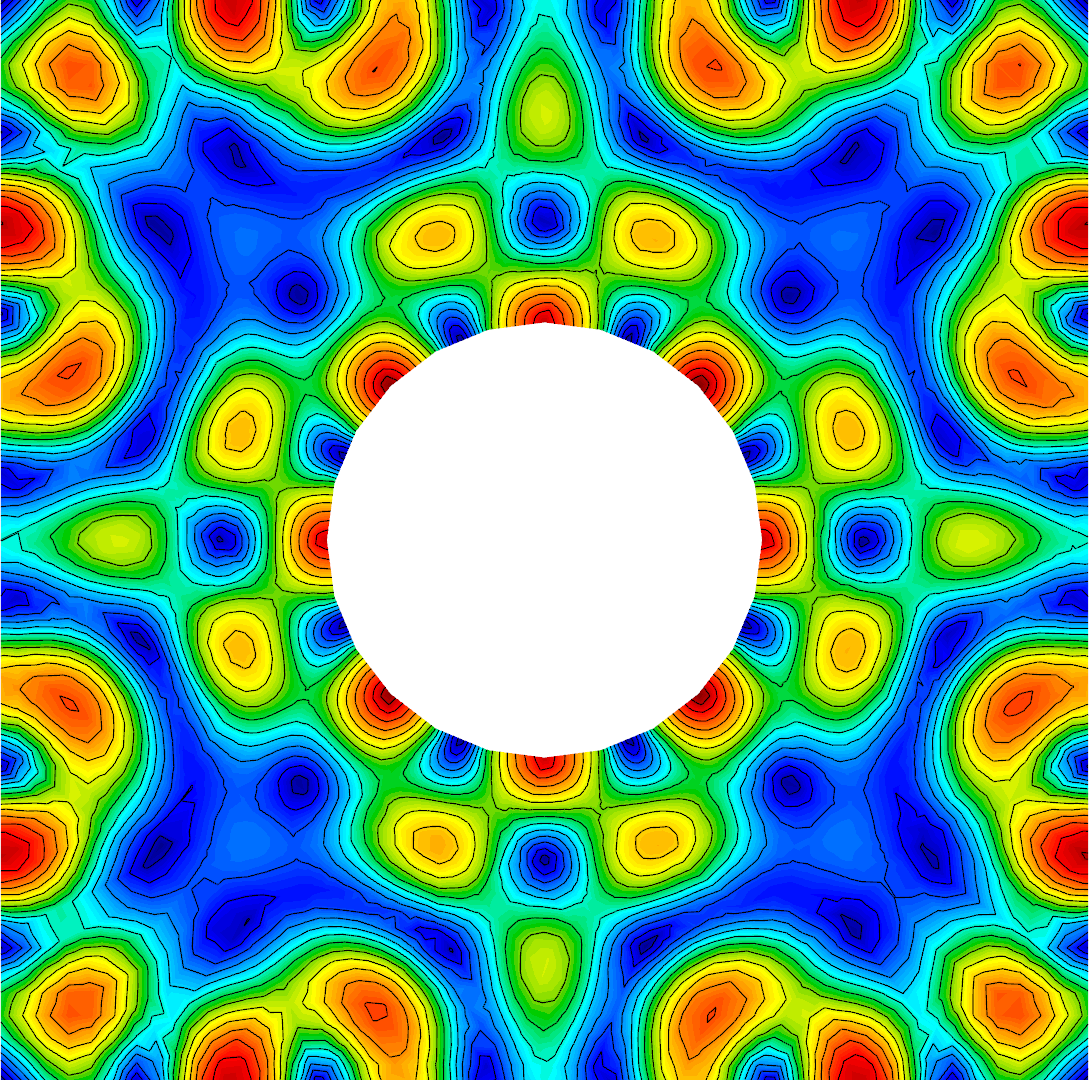} \hspace{0.2em}
\includegraphics[width=0.1\textwidth]{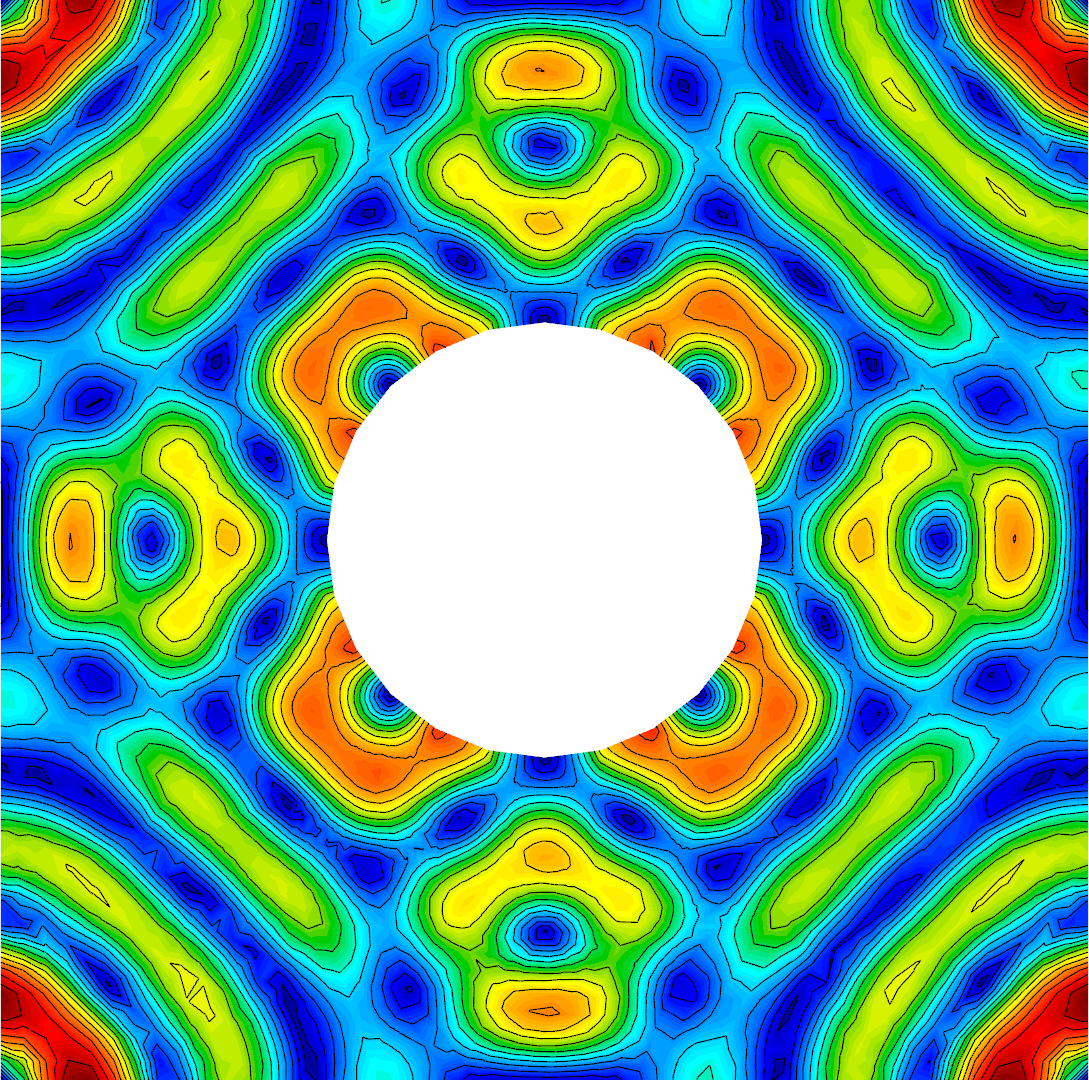} \hspace{0.2em}
\includegraphics[width=0.1\textwidth]{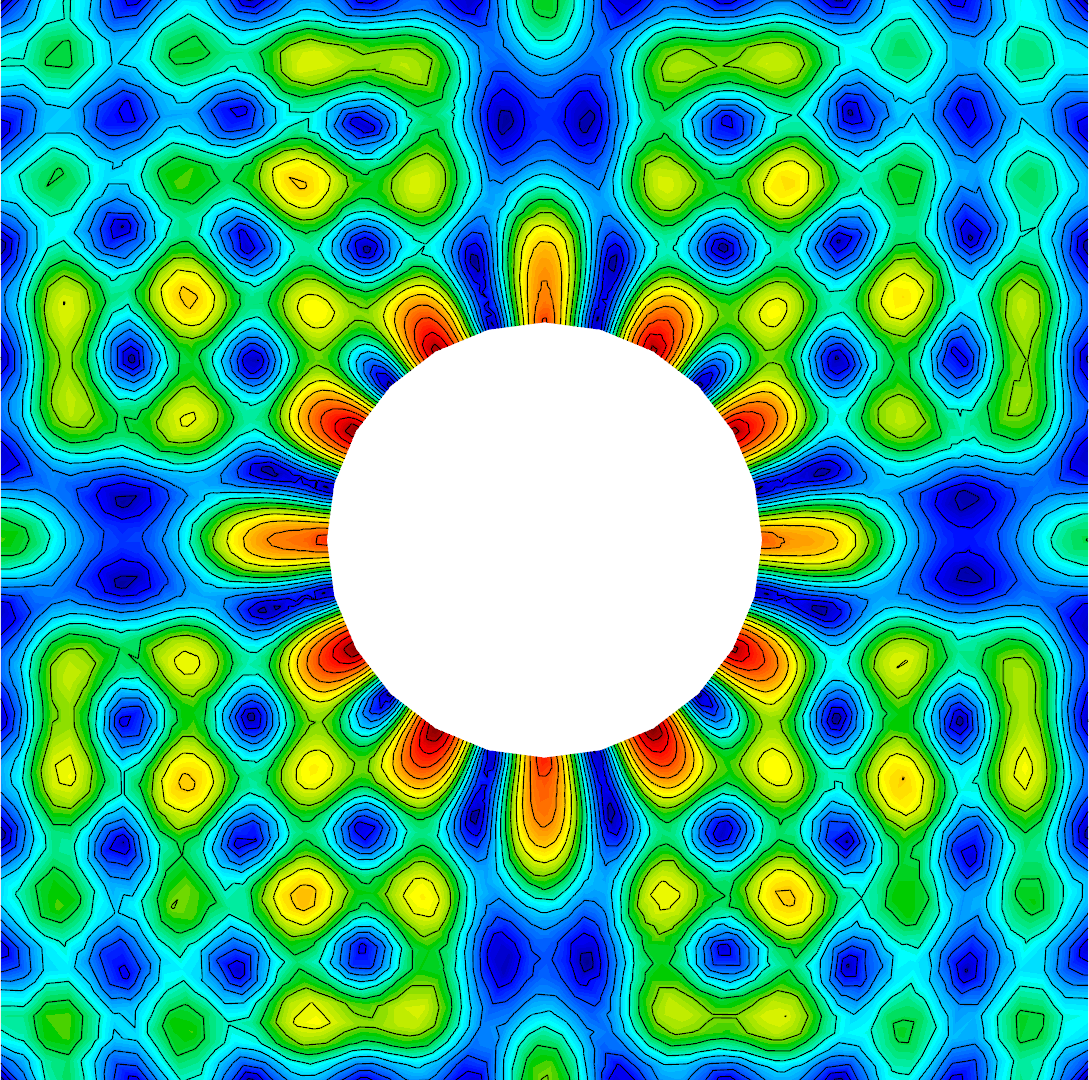} \hspace{0.2em}
\includegraphics[width=0.1\textwidth]{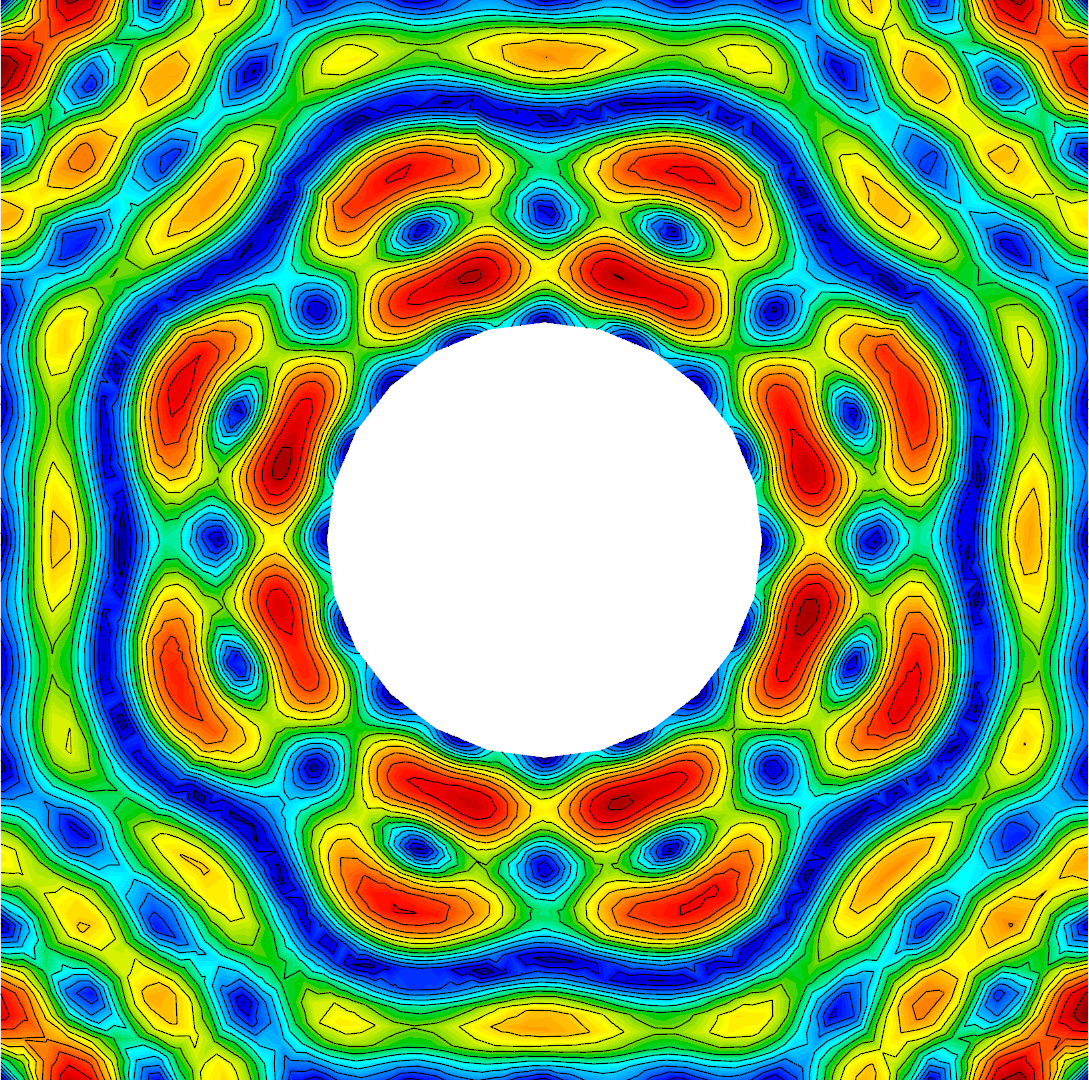} \\[1em]
\includegraphics[width=0.1\textwidth]{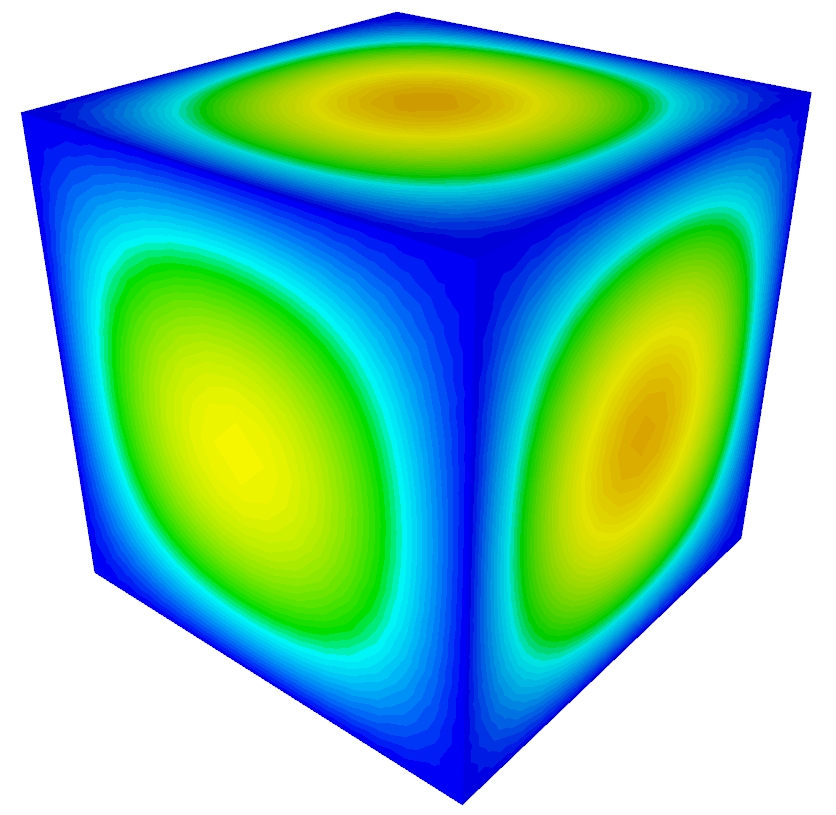} \hspace{0.2em}
\includegraphics[width=0.1\textwidth]{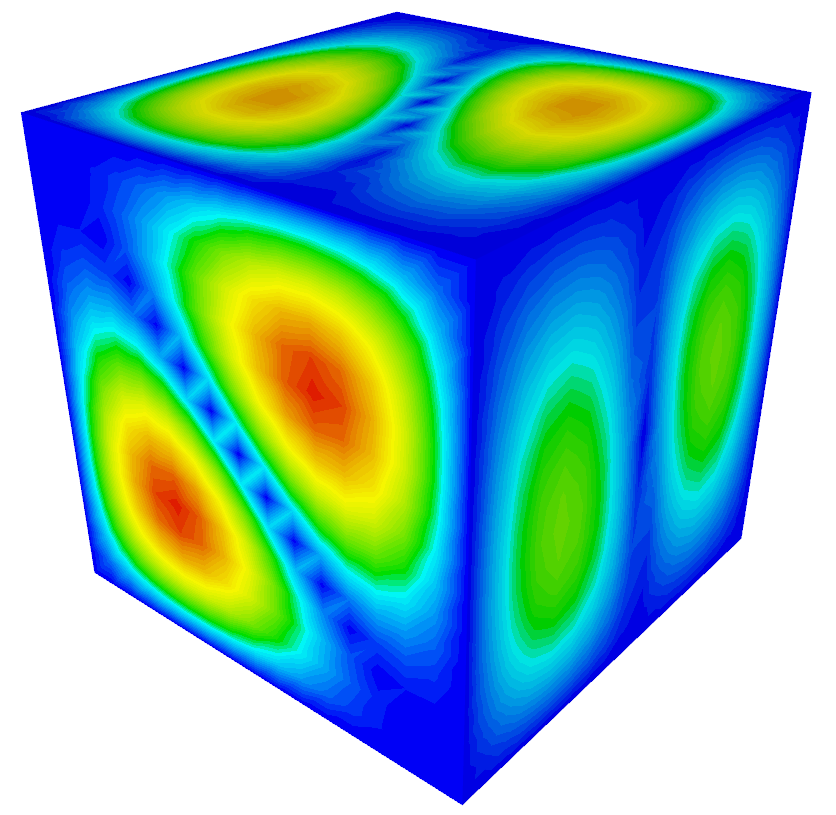} \hspace{0.2em}
\includegraphics[width=0.1\textwidth]{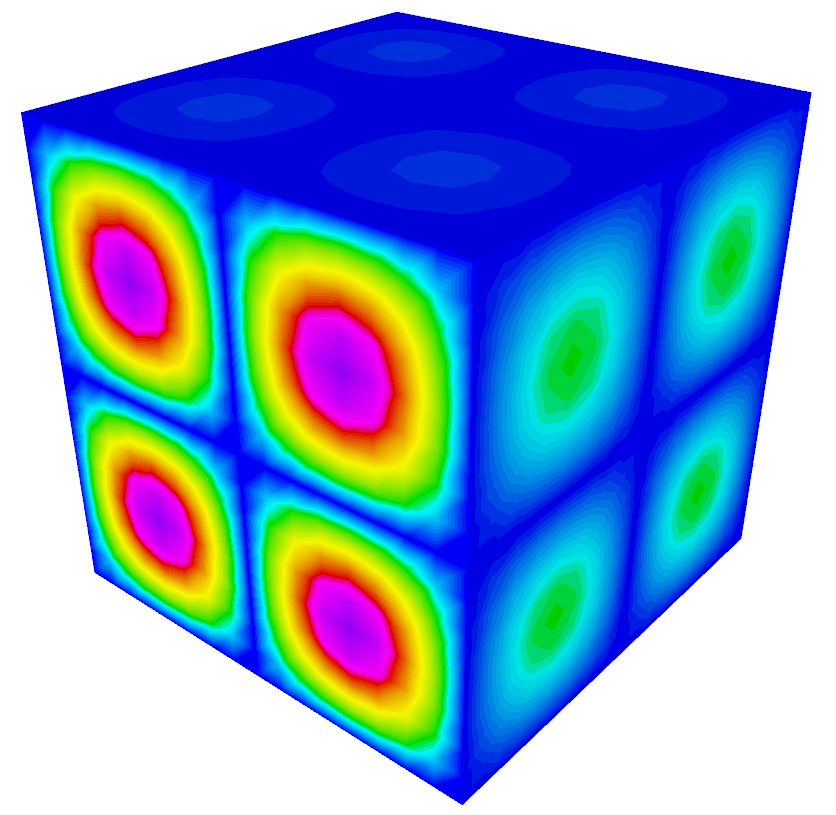} \hspace{0.2em}
\includegraphics[width=0.1\textwidth]{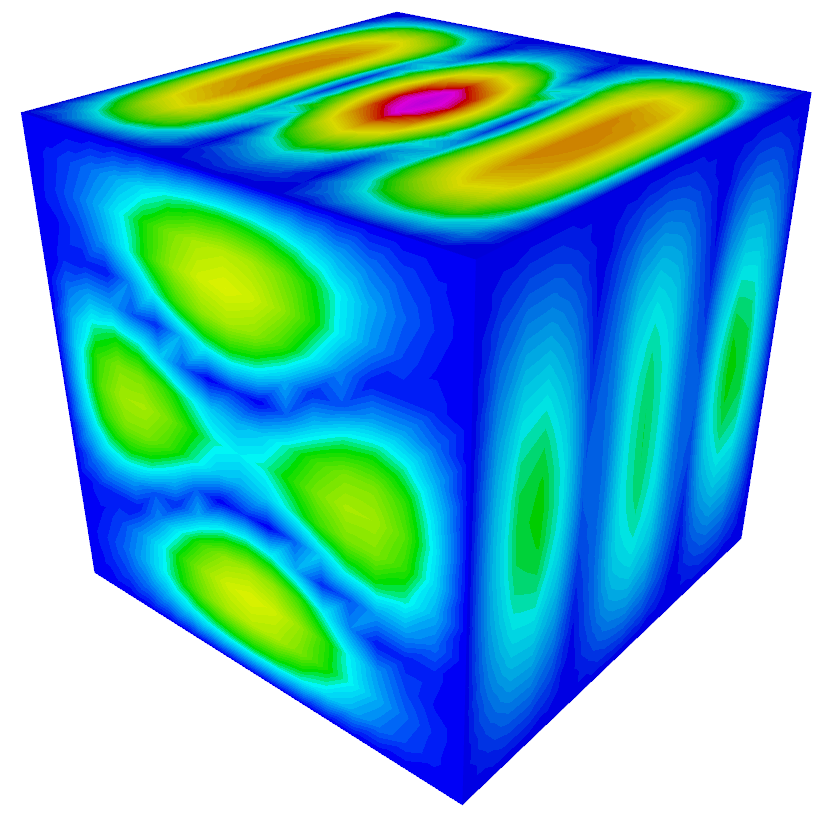} \hspace{0.2em}
\includegraphics[width=0.1\textwidth]{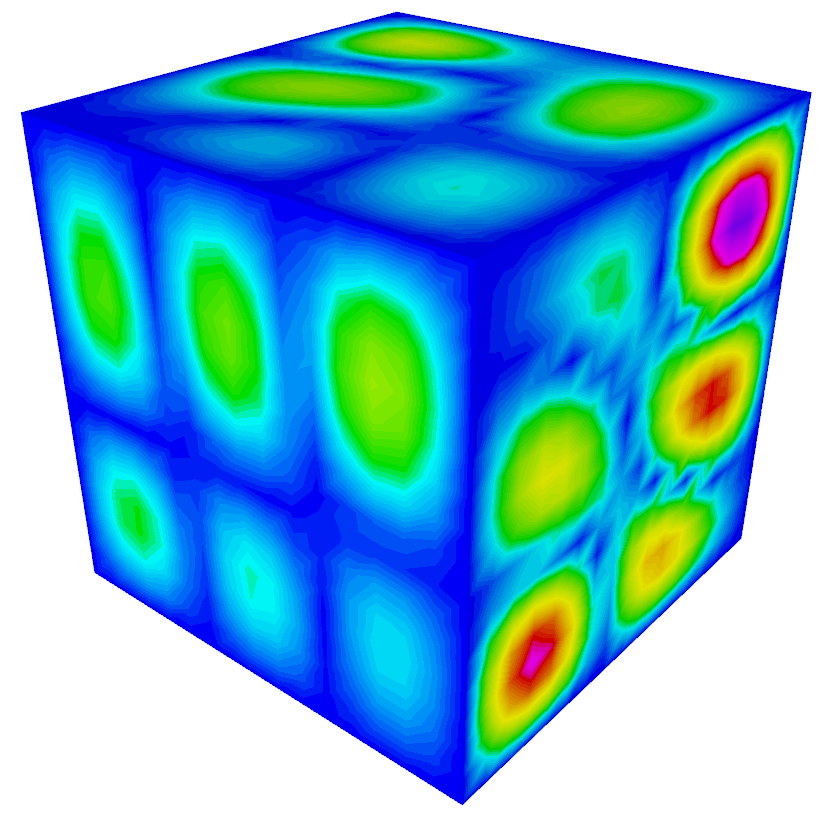} \hspace{0.2em}
\includegraphics[width=0.1\textwidth]{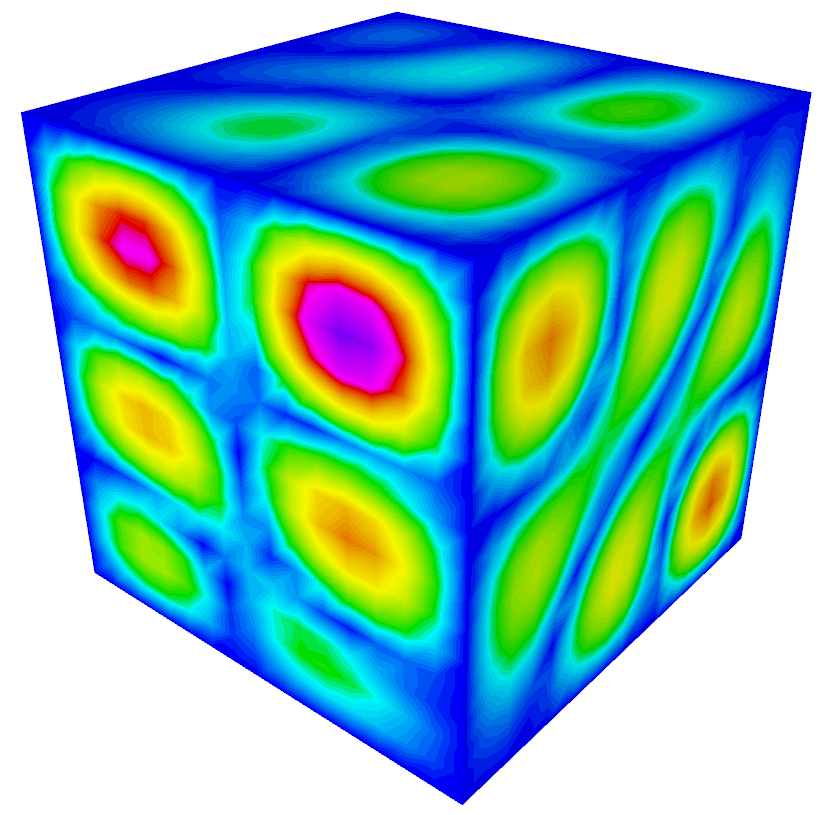} \hspace{0.2em}
\includegraphics[width=0.1\textwidth]{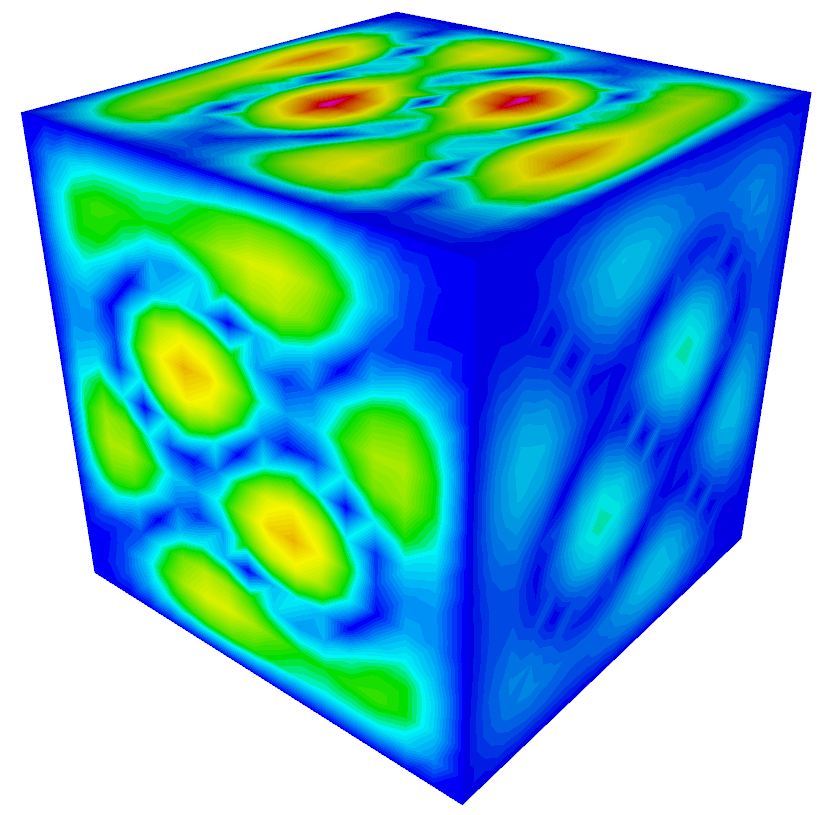} \hspace{0.2em}
\includegraphics[width=0.1\textwidth]{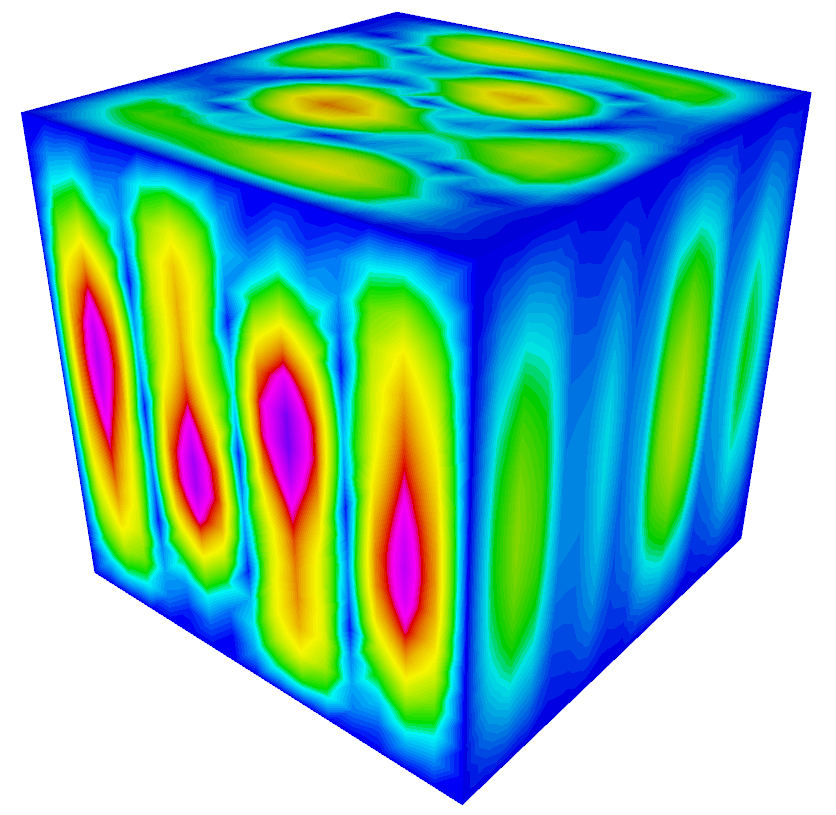} 
\end{figure}

We consider the  Maxwell electromagnetic eigenvalue problem with homogeneous
 Dirichlet boundary conditions:
\begin{align} \label{eq:maxwell}
\nabla \times \nabla \times \vec{E} &= \lambda \vec{E}, \quad \lambda = \omega^2/c^2, \quad \mbox{in} \quad \Omega, \notag \\
\nabla \cdot \vec{E} &= 0, \quad \mbox{in} \quad \Omega, \notag \\
\vec{E} \times \vec{n} &= 0,  \quad \mbox{on} \quad \partial\Omega ,
\end{align}
where $\vec{E}$ is  the electric field, $\vec{n}$ is  the unit outward
normal   to   the   boundary   $\partial\Omega$,   $\omega$   is   the
eigenfrequency  of the  electromagnetic oscillations,  and $c$  is the
light  velocity.   The  curl-curl  operator was  discretized  using  a
N\'{e}d\'{e}lec finite element space of the  second order in 2D or 3D.
Our  goal is  to  find  the electromagnetic  eigenmodes  that are  the
non-zero    solutions   of    \eqref{eq:maxwell}.     As   shown    in
\cite{arbenz2001solving},  the discretized  problem \eqref{eq:maxwell}
is equivalent to the matrix  eigenvalue problem $Ax=\lambda Bx$, where
$A$ is  the stiffness matrix  corresponding to the  curl-curl operator
and $B$ is the mass matrix. It is well-known that the matrix $A$ has a
high-dimensional    zero    eigenspace    spanned    by    discretized
gradients. Therefore, the desired positive eigenvalues are deep inside
the  spectrum.   In   Figure~\ref{fig:maxwellspectrum},  we  show  the
spectrum of  $B\inv A$  for a  2-D problem  on a  square disc  and the
spectrum for  a 3-D  problem on a  cube.  The 2-D  problem is  of size
$3600$ (after  eliminating the  degrees of  freedom on  the boundary),
where  about $32\%$  of the  eigenvalues  ($1137$ out  of $3600$)  are
clustered around  zero and the  smallest eigenvalue that is  away from
zero is  $6.71$.  The 3-D  problem is of  size $10800$ which  also has
about $32\%$ of the eigenvalues ($3375$  out of $10800$) close to zero
and the smallest eigenvalue that is  away from zero is $19.80$.  A few
of  the  non-zero  eigenmodes  of both  problems  are  illustrated  in
Figure~\ref{fig:eigmodes}.

In cases where the knowledge of the  kernel space of $A$ is known from
the                         discretization,                        say
$\mbox{Ker}(A)=\mbox{span}\left\{y_1,\ldots,y_m\right\}$,         this
information can be  used to design more  efficient specialized Maxwell
eigensolvers. For instance, in the Auxiliary Maxwell Eigensolver (AME)
\cite{kolevAME},  the $B$-orthogonal  projector  \eq{eq:ame} P  = I  -
Y(Y\trans  B  Y)\inv Y\trans  B,  \en  where $Y=[y_1,\ldots,y_m]$,  is
applied  within  LOBPCG to  force  the  iterations  to remain  in  the
subspace  $\mbox{Ran}(Y)^{\perp_B}$.   The  same projection  was  also
exploited  in \cite{arbenz2001solving}  as the  preconditioner of  the
implicitly  restarted   Lanczos  algorithm  and   the  Jacobi-Davidson
algorithm.  However, we remark here  that these solvers are not purely
algebraic solvers since additional information from the discretization
is assumed.  On the other hand, with EVSL's filtering mechanism, it is
straightforward to skip  the zero eigenvalues and jump  to the desired
part of the spectrum.
In Table \ref{tab:maxprobs}, we list the 8 Maxwell eigenvalue problems
considered, for which we computed all the eigenvalues 
in the given interval ${[\xi,\eta]}$ and the corresponding eigenvectors.

\begin{table}[tbh]
\caption{Maxwell eigenvalue problems \label{tab:maxprobs}}
\centering
\def\arraystretch{1.1}
\begin{tabular}{r|rrrrr}
 Problem & \multicolumn{1}{c}{n} & \multicolumn{1}{c}{nnz} & \multicolumn{1}{c}{$(a,  b)$} & \multicolumn{1}{c}{$(\xi,  \eta)$} & \multicolumn{1}{c}{$\nu_{[\xi,\eta]}$} \tabularnewline
 \hline
 Max2D-1 & 14,592 & 72,200 & $(0,  9.64 \times 10^{5})$ & $(6,  1200)$ & $96$ \tabularnewline
 Max2D-2 & 59,520 & 292,216 & $(0, 4.00 \times 10^{6})$ & $(6, 1200)$  & $96$
 \tabularnewline 
 Max2D-3 & 235,776 & 1,175,816 & $(0, 1.62 \times 10^{7})$ & $(6, 1200)$ & $96$ \tabularnewline
 Max2D-4 & 944,640 & 4,717,064 & $(0, 6.50 \times 10^{7})$ & $(6, 1200)$ & $96$ \tabularnewline 
 \hline
 Max3D-1 & 10,800 & 320,952 & $(0,  8.96 \times 10^{3})$ & $(19.5,  250)$ & $115$ \tabularnewline
 Max3D-2 & 92,256 & 2,893,800 & $(0, 3.66 \times 10^{4})$ & $(19.5,  250)$ & $121$ \tabularnewline 
 Max3D-3 & 762,048 & 24,526,920 & $(0, 1.47 \times 10^{5})$ & $(19.5,  250)$ & $121$ \tabularnewline
 Max3D-4 & 2,599,200 & 84,363,432 & $(0, 3.32 \times 10^{5})$ & $(19.5,  250)$ &$121$ \tabularnewline 
\end{tabular}
\end{table}

We first present the results with \texttt{ChebLanNr} in Table \ref{tab:maxpol}. 
For the solves with $B$, which is a very well-conditioned mass matrix,
we considered using the sparse direct solver 
Pardiso and the Chebyshev polynomial expansion by \texttt{lsPol}.
These two methods  are indicated as `d' and `c'
in the column  labeled `sv' of Table \ref{tab:maxpol} respectively.
For each problem we present the result with the method for the 
solves with $B$ that yielded
the better total iteration time. We found that for larger problems, in
both 2-D and 3-D, using the Chebyshev iterations yielded
better iteration times, while for the 3 smaller problems using the 
direct solver was more efficient.
In the column `niter', we give
the numbers of the iterations required by the filtered Lanczos algorithm,  which remained roughly the same as the problem sizes increase and thus,
so did the convergence rates.
In order to keep such constant convergence rates, the required degree of the polynomial filter (shown in column `deg') needed to increase, since the target interval $[\xi,\eta]$ 
was unchanged but the spectrum of $B\inv A$ becomes wider as the problem size increases. 
For the iteration time, the cost of performing the solves with
$B$ (t-sv) 
dominated the total iteration time (t-tot), which
is much more expensive than the cost of computing the 
matrix-vector products with $A$ (t-mv) and the cost of the orthogonalization 
(t-orth). 

\begin{table}[tbh]
\caption{Numerical results of  the polynomial filtered Lanczos algorithm for  Maxwell eigenvalue problems \label{tab:maxpol}}
\centering
\tabcolsep1.9mm
\def\arraystretch{1.1}
\begin{tabular}{c|rrcrrrrrr}
\multirow{2}{*}{Problem} & \multirow{2}{*}{deg} & \multirow{2}{*}{niter}  & \multirow{2}{*}{sv} & \multicolumn{2}{c}{fact} & \multicolumn{4}{c}{iter}\tabularnewline
& & & & \multicolumn{1}{c}{mem} & \multicolumn{1}{c}{time} 
& \multicolumn{1}{c}{t-mv} & \multicolumn{1}{c}{t-sv} & 
\multicolumn{1}{c}{t-orth} & \multicolumn{1}{c}{t-tot}\tabularnewline
\hline 
Max2D-1  & 69 & 350  & d & 5 MB & 0.19 & 0.90 & 31.7 & 0.32 & 35.24 \tabularnewline
Max2D-2  & 141 & 350 & d & 21 MB & 0.24 & 2.75 & 221.0 & 1.89 & 233.63 \tabularnewline
Max2D-3  & 284 & 350 & c & - & - & 22.89 & 616.50 & 7.07 & 946.62 \tabularnewline
Max2D-4  & 570 & 350 & c & -& - & 279.20 & 10238.2 & 45.39 & 10585.81 \tabularnewline
\hline 
Max3D-1  & 15 & 530 & d & 32 MB & 0.24 & 0.34 & 9.62 & 0.71 & 11.42 \tabularnewline
Max3D-2  & 30 & 570 & c & - & - & 6.60 & 185.16 & 8.08 & 203.33 \tabularnewline
Max3D-3  & 61 & 590 & c & - & - & 137.67 & 4907.82 & 78.73 & 5149.92 \tabularnewline
Max3D-4  & 93 & 570 & c & - & - & 774.50 & 27407.41 & 529.49 & 28822.21\tabularnewline
\end{tabular}
\end{table}

Next, we examine the performance of \texttt{RatLanNr} shown in Table~\ref{tab:maxrat}.  
Pardiso was used to solve the linear systems with $A-\sigma_j B$ required by
the rational filter. Compared with \texttt{ChebLanNr},
fewer iterations were required and the total iteration time was also much shorter.
For the largest 2-D problem, the iteration time required by \texttt{RatLanNr} 
was $90$ times faster, while for the largest 3-D problem, this speedup was about $3$.
In the previous sections, we have seen significant CPU time savings 
by using \texttt{ChebLanNr} for the 3-D Laplacians and Hamiltonians
compared with \texttt{RatLanNr}.
However, the results presented in this section 
indicate the opposite. The reason for the  different 
time efficiency of the two types of methods for standard and generalized
eigenvalue problems can be understood from the computations required to
 apply the filter at each step of the Lanczos iterations that are given
in Table~\ref{tab:itercomp}. 
Let $k$ denote the degree of the polynomial filter and $k_j$ denote  the number of the repetition
of pole $j$ of the rational filter.
In our runs, only one pole, repeated twice, was used, so we have $\sum k_j = 2$.
For standard eigenvalue problems, each step of \texttt{ChebLanNr} 
requires $k$ matrix-vector products with $A$, 
whereas each step of \texttt{RatLanNr} needs to
perform $\sum k_j$ solves with $A-\sigma_j B$. 
For  large 3-D problems, it turned out that it was
more expensive to perform the solves than the matrix-vector products, so
\texttt{ChebLanNr} was usually more efficient.
However, for the generalized problems, 
there are additional $k$ solves with $B$ for each step of \texttt{ChebLanNr}, which
makes the application of the polynomial filter significantly much more expensive.
In contrast, \texttt{RatLanNr}  requires  $\sum k_j$ 
matrix-vector products with $B$ additionally, 
the cost of which is usually negligible.
Consequently, \texttt{RatLanNr} was often found to be much more efficient for 
generalized eigenvalue problems than \texttt{ChebLanNr}.

\begin{table}[tbh]
\caption{Numerical results of  the rational filtered Lanczos algorithm for  Maxwell eigenvalue problems \label{tab:maxrat}}
\centering
\def\arraystretch{1.1}
\begin{tabular}{c|rcrrrrrr}
\multirow{2}{*}{Problem} & \multirow{2}{*}{niter} & \multirow{2}{*}{sv} 
& \multicolumn{2}{c}{fact} & \multicolumn{4}{c}{iter}\tabularnewline
& & & \multicolumn{1}{c}{mem} & \multicolumn{1}{c}{time} & 
\multicolumn{1}{c}{t-mv} & \multicolumn{1}{c}{t-sv} & 
\multicolumn{1}{c}{t-orth} & \multicolumn{1}{c}{t-tot} \tabularnewline
\hline 
Max2D-1 & 190   & d & 9 MB & 0.10 & .04 & .95 & 0.07 & 1.15 \tabularnewline
Max2D-2 & 190   & d & 39 MB  & 0.28 & .07 & 4.36 & 0.74 & 5.57 \tabularnewline
Max2D-3 & 190   & d & 168 MB & 1.02 & 0.25 & 24.34 & 2.98 & 29.00 \tabularnewline
Max2D-4 & 190   & d & 716 MB & 4.31 & 1.39 & 100.55 & 10.52 & 117.11 \tabularnewline
\hline 
Max3D-1 & 290   & d & 95 MB & 0.23 & 0.06 & 4.95 & 0.18 & 5.39 \tabularnewline
Max3D-2 & 290   & d & 923 MB & 2.45 & 0.46 & 94.65 & 2.28 & 98.20\tabularnewline
Max3D-3 & 290   & d & 13 GB & 64.53 & 5.31 & 1599.72 & 20.69 & 1631.30\tabularnewline
Max3D-4 & 290 & d & 68 GB & 649.21 & 25.21 & 9419.12 & 127.74 & 9609.42 \tabularnewline
\end{tabular}
\end{table}

On the other hand, the memory requirement of \texttt{RatLanNr} is much
higher than \texttt{ChebLanNr}, as shown  in the columns labeled `mem'
of Tables~\ref{tab:maxpol}-\ref{tab:maxrat},  due to the  large memory
consumption  of   the  factorization  of  $A-\sigma_j   B$.   For  the
polynomial filtering, we can avoid any factorization by performing the
solves  with  $B$  with  iterative   methods  such  as  the  Chebyshev
polynomial  iterations.  For   rational  filtering, EVSL  currently
relies on direct methods to  solve the linear systems with the matrices
 $A-\sigma_j B$,  which are highly  indefinite. This is because it  is  challenging 
to  find efficient iterative methods to solve linear systems with  such matrices.


\section{Recommendations, outlook,  and closing remarks} 
\paragraph{Recommendations: Rational  vs. polynomial filtering}
  Based  on our experiments,  we  found that 
  when combined with the Lanczos algorithm, 
  rational filtering tends to  be more efficient
  than  polynomial  filtering,  for  2-D problems  in both  the
  standard and generalized forms. This is because the factorization of
  the  shifted matrix  $A-\sigma_j B$  is generally  inexpensive.  For
  large  3-D  problems in  the  standard  form, the  factorization  of
  $A-\sigma_j B$ becomes  too expensive and in these  cases using the
  polynomial filtered algorithm tends to be more efficient in terms of
  both  CPU  time consumption and memory usage. This can
  be also true for generalized  problems in the situation where linear
  systems with $B$ are (very) inexpensive to solve.
  Finally, for  the situations where  solving linear systems  with the
  $B$ matrix  is costly,  the rational  filtered algorithms  may become
  more  efficient   in  terms   of  CPU   time  usage,   provided  the
  factorization  of $A-\sigma_j  B$ is  still affordable,  although it
  will usually require much more memory.
  Otherwise,  when  factoring  $A-\sigma_j  B$  becomes  prohibitively
  expensive, then the alternative of the polynomial filtered algorithm
  remains a feasible option as long  as linear systems with the matrix
  $B$  can  be efficiently  solved.  There  are situations  in  finite
  elements discretizations  that lead  to mass  matrices $B$  that are
  very well conditioned  once they are scaled by  their diagonals.  In
  these cases, solving a linear system with $B$
and $B^{1/2}$ is quite inexpensive as it can be carried out with 
a small number of matrix-vector
products with $B$, see~\cite{gdos} for details. 

\paragraph{Outlook: Parallelism and iterative solvers} 
A fully parallel version  of EVSL using MPI is currently 
being developed. It is  being tested on a large 
scale geophysics simulation and will be released soon. 
Obtaining a  scalable parallel  rational filtered
algorithm in  a  parallel  distributed memory  environment 
may be rather challenging because of the limitations of current
parallel direct solvers.  
Polynomial filtered algorithms are more  appealing from this perspective since
for standard eigenvalue  problems, the matrix $A$ is  only involved in
matrix-vector products, which can  be efficiently parallelized. Moreover,
as already mentioned,
for  generalized  eigenvalue  problems in which the $B$ matrix is 
well-conditioned, solving systems with $B$ and $B^{1/2}$ can 
be achieved via the
application of  low degree polynomials in $B$ to the right-hand side.

We have not yet mentioned  the possibility of replacing direct
solvers by iterative ones for rational filtering methods, but this is
in  our future plans. Iterative solvers that will be used in this 
context must be able to handle highly indefinite systems, and
methods derived from  those in \cite{XiYS_ratprec,YxiAl2014}
may be suitable. These methods are highly parallel and  must 
 be tuned to the specificity of rational filtering: complex matrices,
availability of spectral information, multiple  shifts, etc.


\paragraph{Closing remarks}
There were attempts  more than two decades ago to  utilize ideas based
on polynomial  filtering for solving interior  eigenvalue problems but
these  were abandoned  because they  were  generally found  not to  be
competitive  relative  to  the standard  Krylov  subspace  approaches.
Similarly, when the  first idea of using Cauchy  integrals came about,
it was not adopted right away.  With the emergence in  recent
years of very large scale
electronic  structure calculations,  and  other scientific  simulations, 
methods  that combine  filtering and  spectrum slicing have become not
just appealing,  but also  mandatory.  Their main  appeal lies
in the  extra level  of parallelism they  offer as well as in
their intrinsic efficiency due to the reduced orthogonalization costs that
characterize  them. Therefore,  the reborn versions of spectrum slicing 
that have recently emerged constitute  a new paradigm that 
is  likely to gain in importance.

\bibliographystyle{siam}
\bibliography{local.bib,nano.bib}

\end{document}